\documentclass[12pt,reqno]{amsart}

\usepackage[T1]{fontenc}
\usepackage[utf8]{inputenc}
\usepackage{lmodern}

\usepackage[a4paper,margin=1in]{geometry}
\usepackage{setspace}

\usepackage{amsmath,amssymb,amsthm,mathtools,amsfonts,bm}
\numberwithin{equation}{section}

\usepackage{graphicx}
\usepackage{subcaption}
\usepackage{booktabs}
\usepackage{multirow}
\usepackage{tabularx}
\usepackage{colortbl}
\usepackage{float}
\usepackage{comment}

\usepackage{tikz}
\usetikzlibrary{arrows.meta,calc,positioning}
\usepackage{pgfplots}
\pgfplotsset{compat=1.18}

\usepackage{enumitem}

\usepackage{lscape}

\usepackage{hyperref}
\hypersetup{
  colorlinks=true,
  linkcolor=blue!50!black,
  citecolor=blue!50!black,
  urlcolor=blue!50!black,
  pdfauthor={},
  pdftitle={}
}
\usepackage[nameinlink,capitalise,noabbrev]{cleveref}

\theoremstyle{plain}
\newtheorem{theorem}{Theorem}[section]
\newtheorem{proposition}[theorem]{Proposition}
\newtheorem{lemma}[theorem]{Lemma}
\newtheorem{corollary}[theorem]{Corollary}

\theoremstyle{definition}

\newtheorem{definition}[theorem]{Definition}
\newtheorem{example}[theorem]{Example}

\theoremstyle{remark}
\newtheorem{remark}[theorem]{Remark}

\DeclareMathOperator*{\Med}{Med}

\newcommand{\R}{\mathbb{R}}
\newcommand{\Rd}{\mathbb{R}^d}

\newcommand{\E}{\mathbb{E}}

\newcommand{\PP}{\mathbf{P}}
\newcommand{\EE}{\mathbf{E}}
\newcommand{\1}{\mathbf{1}}
\newcommand{\Aum}{\mathbb{E}}

\DeclareMathOperator*{\med}{med}

\DeclareMathOperator{\Sel}{\mathbf{Sel}}
\title[Bounds for Restricted Selections]{Bounds for Restricted Selections of Random Sets}
\author[A.~Beresteanu]{Arie Beresteanu}
\address{Department of Economics, University of Pittsburgh, Pittsburgh, PA 15260, USA}
\email{arie@pitt.edu}

\author[B.~M.~Ramezanzadeh]{Behrooz Moosavi Ramezanzadeh}
\address{Department of Economics, University of Pittsburgh, Pittsburgh, PA 15260, USA}
\email{behroozmoosavi@pitt.edu}
\date{\today}

\begin{document}
\begin{abstract}
We study constrained selection sets of random closed sets defined on a
non-atomic probability space. Given a random interval $Y=[y_L,y_U]$
and scalar constraints on the expectation or the median of admissible
selections, we characterize the restricted selection set and establish
sharp bounds on the attainable ranges of means, medians, and event
probabilities. In particular, we give conditions under which every
value in the Aumann expectation range is realized as the mean of a
measurable selection, and we obtain explicit formulas for the extremal
expectations under median and higher-moment restrictions via
rearrangement and convex-duality arguments. We further show that the
selection set of any random compact convex set in $\R^d$ can be
approximated in $L^1$ by selection sets of disjoint unions of random
cubes, each of which decomposes coordinate-wise into one-dimensional
interval selection problems. This gives us an approximation-based reduction of constrained selection problems for random compact convex sets in $\R^d$.
\end{abstract}
\maketitle
\vspace{-15mm}
\tableofcontents

\newpage

\section{Introduction}\label{sec:intro}

Let $(\Omega,\mathcal A,\PP)$ be a complete non-atomic probability space and let
$Y:\Omega\to\mathcal{K}_c(\R^d)$ be a random compact convex set. A measurable selection of $Y$ is a random vector $y:\Omega\to\R^d$ such that
$y(\omega)\in Y(\omega)$ for $\PP$–a.s. The family of all measurable selections is denoted $\Sel(Y)$, and the Aumann expectation of $Y$ is the closed set
\begin{equation}\label{equ:AumannaExpectation}
    \E[Y]:=\overline{\{\E[y]:y\in\Sel^1(Y)\}},
\end{equation}
where $\Sel^1(Y)$ denotes the subfamily of absolute integrable selections; see
\cite{Aumann1965,Artstein1974,Molchanov2005}.

Definition (\ref{equ:AumannaExpectation}) sets the expectation of a random set to be the set of all expectations of its selections. Similar concepts for the median of a random set or for event probabilities can be defined as collections of these operators for all possible selections of the random set (see Subsection \ref{subsec:notation}). In this paper, we study \emph{restricted selection sets} --- subsets of $\Sel(Y)$ defined by imposing scalar constraints on the selections. Concretely, we consider sets of the form
\[
\Sel(Y\mid \kappa)
:=\{y\in\Sel(Y):\E[y]=\kappa\},\qquad
\Sel(Y\mid m):=\{y\in\Sel(Y): m\in\Med(y)\}.
\]
and more generally, selection sets restricted by moment or quantile conditions. Our aim is to understand (i) when such restricted selection sets are nonempty, (ii) how the constraint narrows the set of selections relative to $\Sel(Y)$, and (iii) which values of basic
scalar functionals,
\[
\E[y],\qquad med(y),\qquad \PP(y\in A),
\]
can be attained by selections $y$ satisfying a given constraint.

We develop the core theory in the one-dimensional setting, where $Y=[y_L,y_U]$ is a random interval with $y_L\le y_U$ $\PP$-a.s. In this case the Aumann expectation has the simple representation $\E[Y]=[\E(y_L),\E(y_U)]$, so that for every $y\in\Sel^1(Y)$,
\begin{equation}\label{equ:expectationBounds}
    \E(y_L) \leq \E[y] \leq \E(y_U).
\end{equation}
Similarly, since each $y\in\Sel(Y)$ stochastically dominates $y_L$ and is stochastically dominated by $y_U$,
\begin{equation}\label{equ:medianBounds}
    med(y_L) \leq med(y) \leq med(y_U).
\end{equation}
The central question is: how do the attainable ranges in \eqref{equ:expectationBounds} and \eqref{equ:medianBounds} change when one restricts attention to a nonempty subset $\mathcal{S}\subset\Sel(Y)$ defined by a scalar constraint?

The non-atomicity of $(\Omega,\mathcal A,\PP)$ allows us to use measurable selection and
rearrangement arguments in a sharp way. At a technical level, several of our proofs
combine classical properties of the Aumann expectation with Hardy–Littlewood type
inequalities and quantile representations. This yields explicit and, in many cases,
closed-form descriptions of the extremal values of expectations, medians, event
probabilities, and higher moments over constrained selection sets.

We then show that the one-dimensional theory is sufficient for the general case: the selection set of any random compact convex set in $\R^d$ can be approximated in $L^1$ by selection sets of disjoint unions of random cubes, each of which decomposes coordinate-wise into a product of one-dimensional interval selection sets. This reduces every constrained selection problem for random convex sets in $\R^d$ to the interval case.

\subsection{Related Literature}

Our work is rooted in the classical theory of random closed sets and selection
expectation. The notion of the Aumann integral for set-valued mappings was introduced
in \cite{Aumann1965} and developed further in \cite{Artstein1974}, where convexity and
support-function representations play a central role. In the one-dimensional setting
considered here, the random sets are intervals $Y=[y_L,y_U]$, and the identity
$\E[Y]=[\E(y_L),\E(y_U)]$ is a special case of these general results. For the basic
theory of random closed sets, measurable selections, and the measurable graph theorem, we
refer to \cite{KuratowskiRyll1965,CastaingValadier1977,Molchanov2005}.

Rearrangement and quantile methods appear throughout analysis and probability in the
study of extremal inequalities. Our arguments for quantile attainability and
median-restricted expectations are in the spirit of classical
Hardy-Littlewood-P\'{o}lya inequalities, as presented for instance in
\cite{LiebLoss2001,HardyLittlewoodPolya1952}. We employ these tools in a simple random
interval setting, where they yield explicit formulas for the effect of median
constraints on the expectation and for extremal event probabilities under a mean
restriction.

Random sets and Aumann expectations have been used in a variety of applied contexts,
including mathematical economics, finance, and statistics; see, for example, the
monograph \cite{Molchanov2005} and the references therein. In particular, models with
set-valued or interval-valued observations naturally lead to constrained selection
problems of the type studied here.
The present paper is primarily structural: we develop the core results in a one-dimensional
probability-theoretic setting and then show how they extend to random compact convex sets in $\R^d$ via a dyadic approximation of the selection set.

\subsection{Our Contribution}

We develop a detailed theory of \emph{restricted selection sets} for
random closed sets on a non-atomic probability space. Starting from the benchmark Aumann
expectation range and median interval for the full selection set $\Sel(Y)$, we first
show that these bounds are sharp: every point of
$\E[Y]=[\E(y_L),\E(y_U)]$ is realized as $\E[y]$ for some $y\in\Sel(Y)$
(Proposition~\ref{prop:nonemptySelectionSets}). More generally, for each
$\alpha\in(0,1)$, any value between the capacity and containment quantiles of $Y$ can be
attained as the $\alpha$–quantile of a selection (Proposition~\ref{prop:quantile_selection}).

We then analyze how scalar constraints narrow the selection set and sharpen these bounds. In the median-restricted case, we obtain sharp bounds on the attainable means: under mild feasibility conditions, the set
\[
\{\E[y]:y\in\Sel(Y\mid m)\}
\]
is a compact interval $[\E_{\min}(m),\E_{\max}(m)]$ whose endpoints admit an explicit
representation in terms of quantile integrals of the ``gaps'' $y_U-m$ and $m-y_L$ on the
contact set $\{y_L\le m\le y_U\}$ (Proposition~\ref{prop:MedianRestrictedAumann}). The
proof relies on two auxiliary rearrangement lemmas, including a quantitative
quantile–area identity.

We also examine the converse: does a mean constraint $\E[y]=\kappa$ narrow the set
of attainable medians? A two-state construction
(Example~\ref{ex:no-shrink-med}) shows that the answer is negative in general:
every point in the benchmark median interval
\[
[\med_-(y_L),\,\med_+(y_U)]
\]
can remain attainable as a median under the mean constraint. For event probabilities $\PP(y\in A)$, we refine the classical random-set bounds in Proposition~\ref{prop:RS-unrestricted} by imposing a mean constraint and deriving a
dual (Lagrangian) representation (Theorem~\ref{thm:dual}). This leads to an explicit threshold selection rule (Theorem~\ref{thm:threshold-calibration}) and a closed-form formula for the extremal probabilities $\PP(y\in A)$ over $\Sel(Y\mid\kappa)$ (Corollary~\ref{cor:IdProb}), while preserving consistency with the unrestricted capacity/containment bounds (Proposition~\ref{prop:extremes}).

We further extend the framework in three directions. First, we consider higher-moment constraints and obtain dual representations for the extremal means under an $r$th-moment restriction (Theorem~\ref{thm:r-moment-dual}). Second, we treat quantile constraints at general levels $\alpha\neq\tfrac12$ and show that the induced set of attainable means is always a convex subset of the Aumann range (Proposition~\ref{prop:alpha-quantile-mean}). Finally, we show that the selection set of any random compact convex set $Y\subset\R^d$ can be approximated in $L^1$ by selection sets of disjoint unions of random cubes, giving us an approximation-based reduction of constrained selection problems to finite families of one-dimensional interval-selection problems (Theorem~\ref{thm:reduction}).
\subsection{Notation and preliminaries}\label{subsec:notation}

We work on a fixed, complete, non-atomic probability space
\[
(\Omega,\mathcal{A},\PP).
\]
All random objects are defined on this space. Expectations are denoted by
$\E[\cdot]$, and indicator functions of events $B\in\mathcal{A}$ by
$\mathbf{1}_{B}$. For $x\in\R$ we write $x_+ := \max\{x,0\}$.

\medskip

\noindent\textbf{Spaces and sets.}
We denote by $\R$ the real line and, for $d\in\mathbb{N}$, by $\R^d$ the
$d$–dimensional Euclidean space equipped with its usual topology and Euclidean
norm $\|\cdot\|$. Let $\mathcal{K}(\R)$ be the family of all nonempty closed
subsets of $\R$. Our random sets always take values in $\mathcal{K}(\R)$.
For general background on random closed sets we refer to \cite{Molchanov2005}.

\medskip

\noindent\textbf{Random intervals.}
A \emph{random closed set} in $\R$ is a measurable map
\[
Y:(\Omega,\mathcal{A})\to\big(\mathcal{K}(\R),\mathcal{B}(\mathcal{K}(\R))\big),
\]
where $\mathcal{B}(\mathcal{K}(\R))$ is the Borel $\sigma$–algebra generated by
the Fell (hit–or–miss) topology (see \cite{Molchanov2005}). In this note we
consider the special case of bounded \emph{random intervals}, that is,
\[
Y(\omega) = [y_L(\omega),y_U(\omega)] \in \mathcal{K}(\R),\qquad \omega\in\Omega,
\]
for some real-valued random variables $y_L,y_U:\Omega\to\R$ satisfying
\[
\PP\big(y_L \le y_U\big) = 1.
\]
Unless otherwise stated we assume
\begin{equation}\label{eq:int-moment}
  \E\big(|y_L| + |y_U|\big) < \infty.
\end{equation}

\medskip

\noindent\textbf{Selections and Aumann expectation.}
A (measurable) \emph{selection} of a random closed set $Y$ is a random variable
$y:\Omega\to\R$ such that
\[
y(\omega)\in Y(\omega)\qquad\text{for $\PP$–almost every }\omega.
\]
The existence of measurable selections for closed-valued measurable
multifunctions is classical, see e.g.\ \cite{KuratowskiRyll1965,CastaingValadier1977}.
We denote the family of all measurable selections of $Y$ by
\[
\Sel(Y)
:=
\big\{y:\Omega\to\R\ \text{measurable} : y(\omega)\in Y(\omega)
\ \text{$\PP$–a.s.}\big\},
\]
and the subfamily of integrable selections by
\[
\Sel^1(Y)
:=
\big\{y\in\Sel(Y):\E|y|<\infty\big\}.
\]
Under \eqref{eq:int-moment}, every selection of a random interval $Y=[y_L,y_U]$
is integrable, so $\Sel(Y)=\Sel^1(Y)$ in this case.

Following \cite{Aumann1965}, the \emph{Aumann expectation} (selection expectation)
of a random closed set $Y$ with at least one integrable selection is defined as
the closed set
\[
\E[Y]
:=
\overline{\big\{\E[y]:y\in\Sel^1(Y)\big\}} \subset \R,
\]
where the closure is taken in the usual topology on $\R$; see also
\cite{Artstein1974,Molchanov2005}. For a random interval $Y=[y_L,y_U]$
satisfying \eqref{eq:int-moment}, one has the simple representation
\begin{equation}\label{eq:Aumann-interval}
  \E[Y] = \big[\E(y_L),\,\E(y_U)\big],
\end{equation}
which is a special case of the general convexity and support–function properties
of Aumann integrals.

\medskip

\noindent\textbf{Capacity and containment functionals.}
For a random closed set $Y$ and a Borel set $B\subseteq\R$, the \emph{capacity}
(hitting) and \emph{containment} functionals are defined by
\[
T_Y(B) := \PP\big(Y\cap B\neq\emptyset\big),\qquad
C_Y(B) := \PP\big(Y\subseteq B\big),
\]
see \cite{Molchanov2005}. In the one–dimensional setting we will mainly use the
half–line sets $(-\infty,t]$, $t\in\R$, and write
\[
T_Y((-\infty,t]) = \PP\big(Y\cap(-\infty,t]\neq\emptyset\big),\qquad
C_Y((-\infty,t]) = \PP\big(Y\subseteq(-\infty,t]\big).
\]
Their (left-continuous) generalized inverses are given by
\[
T_Y^{-1}(\alpha)
:= \inf\{t\in\R:T_Y((-\infty,t])\ge\alpha\},\qquad
C_Y^{-1}(\alpha)
:= \inf\{t\in\R:C_Y((-\infty,t])\ge\alpha\},
\]
for $\alpha\in(0,1)$. For random intervals $Y=[y_L,y_U]$ one has
\[
T_Y((-\infty,t]) = \PP(y_L\le t),\qquad
C_Y((-\infty,t]) = \PP(y_U\le t),
\]
so that $T_Y^{-1}$ and $C_Y^{-1}$ coincide with the quantile functions of
$y_L$ and $y_U$, respectively.

\medskip

\noindent\textbf{Distribution functions, quantiles, and medians.}
If $z$ is a real-valued random variable with distribution function
\[
F_z(t) := \PP(z\le t),\qquad t\in\R,
\]
its (left-continuous) quantile function is
\[
F_z^{-1}(\alpha)
:= \inf\{t\in\R:F_z(t)\ge\alpha\},\qquad \alpha\in(0,1).
\]
A (possibly non-unique) median of $z$ is any $m\in\R$ such that
\[
\PP(z\le m)\ge\tfrac12,
\qquad
\PP(z\ge m)\ge\tfrac12.
\]
We write $\Med(z)$ for the set of all such medians. In the scalar case,
$\Med(z)$ is always a nonempty closed interval (possibly a singleton). We also use
\[
\med_-(z):=\inf \Med(z),
\qquad
\med_+(z):=\sup \Med(z)
\]
for the lower and upper endpoints of the median set.

\medskip

\noindent\textbf{Restricted selection sets.}
The objects of interest in this note are selection sets restricted by simple
scalar constraints. For $\kappa\in\R$ we define the \emph{$\kappa$-mean restricted
selection set} by
\[
\Sel(Y\mid \kappa)
:=
\big\{y\in\Sel(Y):\E[y]=\kappa\big\},
\]
and for $m\in\R$ the \emph{$m$-median restricted selection set} by
\[
\Sel(Y\mid m)
:=
\big\{y\in\Sel(Y): m\in\Med(y)\big\}.
\]

\section{Main Results}\label{sec:main}

In this section, we study restricted selection sets of a random interval
$Y=[y_L,y_U]$ on the non-atomic probability space
$(\Omega,\mathcal{A},\PP)$ introduced in Subsection~\ref{subsec:notation}.
First, we record benchmark ranges for the expectation, median, and event probabilities of measurable selections of $Y$. Then we investigate how these ranges are sharpened by simple scalar restrictions on the mean or median of the selection. Throughout we assume $y_L,y_U$ are measurable real-valued random variables with
$\PP(y_L\le y_U)=1$ and
\[
\E(|y_L|+|y_U|)<\infty,
\]
so that every selection $y\in\Sel(Y)$ is integrable and
$\Sel(Y)=\Sel^1(Y)$.

\subsection{Benchmark bounds}\label{subsec:benchmark}

 Let $y\in\Sel(Y)$ be an arbitrary measurable selection of the random interval $Y=[y_L,y_U]$. The following bounds for the expectation and median of a selection $y$ are immediate.

First, by the Aumann expectation identity \eqref{eq:Aumann-interval},
\[
\E[Y]=\big[\E(y_L),\E(y_U)\big],
\]
so the set of possible expectations of selections equals
\[
\big\{\E[y]:y\in\Sel(Y)\big\}
=
\big[\E(y_L),\E(y_U)\big].
\]
Thus
\begin{equation}\label{equ:benchmarkBounds_E}
 \E[y] \in \Theta_E := \big[\E(y_L),\E(y_U)\big].
\end{equation}

Second, since $y_L\le y\le y_U$ almost surely, every median of $y$ must lie between
the lower endpoint of the median set of $y_L$ and the upper endpoint of the median
set of $y_U$. Thus
\begin{equation}\label{equ:benchmarkBounds_M}
  \Med(y)\subseteq \Theta_M
  :=\big[\med_-(y_L),\,\med_+(y_U)\big].
\end{equation}
where the right-hand side is understood as the smallest interval containing both $\Med(y_L)$ and $\Med(y_U)$.

The bounds in \eqref{equ:benchmarkBounds_E} and \eqref{equ:benchmarkBounds_M}
are sharp in the sense that every point in these intervals is the
expectation (respectively, a median) of some measurable selection $y$ of
$Y$. Sharpness for expectations is a direct consequence of
\eqref{eq:Aumann-interval}; sharpness for medians will follow from the quantile attainability result in subsection~\ref{subsec:quantile}.

\subsection{Restricted selection sets}\label{subsec:restricted-sets}

We now introduce restricted selection sets defined by simple scalar
constraints, which will be the central objects of interest.

\begin{definition}[Mean-restricted selection set]\label{def:meanRestrictedSelections}
For $\kappa\in\R$, the $\kappa$-mean restricted selection set is
\[
\Sel(Y\mid \kappa)
:=
\big\{y\in\Sel(Y):\E[y]=\kappa \big\},
\]
and the set $\Sel^1(Y\mid \kappa)$ is defined analogously, replacing
$\Sel(Y)$ by $\Sel^1(Y)$. Under the integrability condition
$\E(|y_L|+|y_U|)<\infty$, these coincide.
\end{definition}

\begin{proposition}\label{prop:nonemptySelectionSets}
Let $Y=[y_L,y_U]$ be a random interval with measurable endpoints
satisfying $y_L\le y_U$ almost surely and $y_L,y_U\in L^1(\PP)$. Then
for any $\kappa \in [\EE(y_L),\EE(y_U)]$,
\[
\Sel(Y\mid\kappa)\neq\emptyset.
\]
\end{proposition}

\begin{proof}
Let $Y=[y_L,y_U]$ be a random interval with measurable endpoints satisfying
$y_L\le y_U$ a.s.\ and assume $y_L,y_U\in L^1(\PP)$, the space of integrable random variables.
For each $t\in[0,1]$ define
\[
y_t \;\coloneqq\; (1-t)\,y_L + t\,y_U.
\]
Then $y_t$ is measurable and, since $y_L\le y_t\le y_U$ a.s., we have
$y_t\in\Sel(Y)$ (indeed $y_t\in\Sel^1(Y)$ by integrability).
By linearity of expectation,
\[
\EE[y_t] \;=\; (1-t)\,\EE[y_L] + t\,\EE[y_U],
\]
so the map $t\mapsto \EE[y_t]$ is affine with range
$[\EE(y_L),\EE(y_U)]$.
Given any $\kappa\in[\EE(y_L),\EE(y_U)]$, choose
\[
t^\ast \;=\;
\begin{cases}
\dfrac{\kappa-\EE(y_L)}{\EE(y_U)-\EE(y_L)}, & \text{if }\EE(y_U)>\EE(y_L),\\[6pt]
0, & \text{if }\EE(y_U)=\EE(y_L),
\end{cases}
\]
and set $y\coloneqq y_{t^\ast}$.
Then $y\in\Sel(Y)$ and $\EE[y]=\kappa$, hence
$\Sel(Y\mid \kappa)\neq\emptyset$ (and, under $L^1$ endpoints,
$\Sel^1(Y\mid \kappa)\neq\emptyset$).
\end{proof}

\begin{definition}[Median-restricted selection set]\label{def:medianRestrictedSelections}
For $m\in\R$, the $m$-median restricted selection set is
\[
\Sel(Y\mid m)
:=
\big\{y\in\Sel(Y): m\in\Med(y)\big\},
\]
and $\Sel^1(Y\mid m)$ is defined analogously.
\end{definition}

Whenever $m$ lies between the capacity and containment medians of $Y$, the
non-emptiness of $\Sel(Y\mid m)$ follows from the quantile attainability
result in Subsection~\ref{subsec:quantile} with $\alpha=\tfrac12$.

The rest of this section describes the induced ranges of expectations, medians,
and event probabilities under these restrictions.

\subsection{Quantile attainability via selections}\label{subsec:quantile}

We now show that any value between the capacity and containment
quantiles of $Y$ can be realized as a quantile of some measurable
selection. Recall from Subsection~\ref{subsec:notation} that for a
random closed set $Y$ and $t\in\R$,
\[
T_Y((-\infty,t]) := \PP\!\big(Y\cap(-\infty,t]\neq\emptyset\big),\qquad
C_Y((-\infty,t]) := \PP\!\big(Y\subset(-\infty,t]\big),
\]
and their generalized inverses
\[
T_Y^{-1}(\alpha):=\inf\{t: T_Y((-\infty,t])\ge\alpha\},\quad
C_Y^{-1}(\alpha):=\inf\{t: C_Y((-\infty,t])\ge\alpha\},\quad \alpha\in(0,1).
\]

\begin{proposition}[Quantile attainability via selections]\label{prop:quantile_selection}
Let $(\Omega,\mathcal F,\PP)$ be a non-atomic probability space and let $Y$ be a
(measurable) random closed convex set in $\mathbb R$ (a random interval). For $t\in\mathbb R$, define
\[
T_Y((-\infty,t]) := \PP\!\big(Y\cap(-\infty,t]\neq\emptyset\big),\qquad
C_Y((-\infty,t]) := \PP\!\big(Y\subset(-\infty,t]\big),
\]
and let the generalized inverses be
\[
T_Y^{-1}(\alpha):=\inf\{t: T_Y((-\infty,t])\ge\alpha\},\quad
C_Y^{-1}(\alpha):=\inf\{t: C_Y((-\infty,t])\ge\alpha\},\quad \alpha\in(0,1).
\]
If $m\in[T_Y^{-1}(\alpha),\,C_Y^{-1}(\alpha)]$, then there exists a measurable selection
$y$ of $Y$ such that $F_y^{-1}(\alpha)=m$.
\end{proposition}

\begin{proof}
Fix $\alpha\in(0,1)$ and $m\in[T_Y^{-1}(\alpha),\,C_Y^{-1}(\alpha)]$, and write $H_t:=(-\infty,t]$.
Recall that for an interval-valued random set $Y=[y_L,y_U]$ we have
\[
T_Y(H_t)=\PP(y_L\le t),\qquad C_Y(H_t)=\PP(y_U\le t),\qquad t\in\R.
\]

\smallskip\noindent
\textbf{Partition at $m$.}
Partition $\Omega$ at the threshold $m$ as
\[
A^- := \{Y\subset H_m\}=\{y_U\le m\},\qquad
A^+ := \{Y\cap H_m=\emptyset\}=\{y_L>m\},
\]
and
\[
A^0 := \Omega\setminus(A^-\cup A^+).
\]
Then
\[
\PP(A^-)=C_Y(H_m),\qquad \PP(A^-\cup A^0)=T_Y(H_m).
\]
Since $m\ge T_Y^{-1}(\alpha)$, $T_Y(H_m)\ge\alpha$, so
\[
\PP(A^0)=T_Y(H_m)-\PP(A^-)\ \ge\ \alpha-\PP(A^-).
\]
Set
\[
\delta:=\bigl(\alpha-\PP(A^-)\bigr)_+ \in[0,\PP(A^0)].
\]
By non-atomicity there exists a measurable $B\subset A^0$ with $\PP(B)=\delta$.

\smallskip\noindent
\textbf{Extremal measurable selections.}
For each $\omega$, define
\[
Y^{\le m}(\omega):=Y(\omega)\cap(-\infty,m],\qquad
Y^{> m}(\omega):=Y(\omega)\cap(m,\infty).
\]
On $A^-\cup A^0$, the set $Y^{\le m}(\omega)$ is nonempty; on $A^0\cup A^+$, the set
$Y^{>m}(\omega)$ is nonempty (note that on $A^0$ we have $y_U(\omega)>m$, since
$y_U(\omega)\le m$ together with $y_L(\omega)\le m$ would place $\omega$ in $A^-$).
Because $Y$ has a measurable graph and the half-lines are Borel,
$\operatorname{graph}(Y^{\le m})$ and $\operatorname{graph}(Y^{>m})$ are measurable.
By standard measurable selection results
(e.g.\ via a Castaing representation, see
\cite{KuratowskiRyll1965,CastaingValadier1977}) we can take measurable
selections
\[
z^-(\omega)\in Y^{\le m}(\omega)\quad\text{on }A^-\cup A^0,\qquad
z^+(\omega)\in Y^{>m}(\omega)\quad\text{on }A^0\cup A^+.
\]

Since $Y(\omega)=[y_L(\omega),y_U(\omega)]$ is an interval, we have the
following explicit forms:
\[
z^-(\omega)=
\begin{cases}
y_U(\omega), & \omega\in A^-,\\[2pt]
m, & \omega\in A^0,
\end{cases}
\]
(noting that $\sup Y^{\le m}(\omega)=m$ on $A^0$ and $=y_U(\omega)$ on $A^-$),
while on $A^0\cup A^+$, every selection $z^+$ of $Y^{>m}$ satisfies
\begin{equation}\label{eq:z-plus-strict}
z^+(\omega)>m.
\end{equation}
On $A^+$ we take $z^+(\omega)=y_L(\omega)>m$.
On $A^0$ we set $z^+(\omega)\in(m,y_U(\omega)]$, so in particular
$z^+(\omega)>m$.

For notational convenience, we complete $z^-$ and $z^+$ to selections of $Y$ by setting
\[
y^-(\omega):=
\begin{cases}
z^-(\omega), & \omega\in A^-\cup A^0,\\[2pt]
y_L(\omega), & \omega\in A^+,
\end{cases}
\qquad
y^+(\omega):=
\begin{cases}
z^+(\omega), & \omega\in A^0\cup A^+,\\[2pt]
y_U(\omega), & \omega\in A^-.
\end{cases}
\]
Thus $y^-,y^+\in\Sel(Y)$ and
\[
y^-(\omega)\le m\ \text{ on }A^-\cup A^0,\qquad
y^+(\omega)> m\ \text{ on }A^0\cup A^+.
\]

\smallskip\noindent
\textbf{Pasting.}
Define $y:\Omega\to\R$ by
\[
y(\omega)=
\begin{cases}
y^-(\omega), & \omega\in A^-\cup B,\\[2pt]
y^+(\omega), & \omega\in (A^0\setminus B)\cup A^+.
\end{cases}
\]
Then $y(\omega)\in Y(\omega)$ for all $\omega$, so $y$ is a measurable selection.

On $A^-\cup B$ we have $y\le m$, and on $(A^0\setminus B)\cup A^+$ we have $y\ge m$, with
$y>m$ on $A^+$. Hence
\[
\PP(y\le m)
\;\ge\;
\PP(A^-)+\PP(B)
\;=\;
\PP(A^-)+\delta
\;=\;
\max\{\PP(A^-),\,\alpha\}
\;\ge\;\alpha.
\]
Thus $F_y(m):=\PP(y\le m)\ge\alpha$.

\smallskip\noindent
\textbf{Behaviour below $m$.}
Let $t<m$. On $A^+$ and on $A^0\setminus B$ we have $y\ge m>t$, so
$\{y\le t\}$ does not occur there. On $B\subset A^0$ we have
$y(\omega)=y^-(\omega)=z^-(\omega)=m>t$, so $\{y\le t\}$ does not occur on $B$ either.
On $A^-$ we have $y(\omega)=y^-(\omega)=z^-(\omega)=y_U(\omega)$, since
$Y(\omega)\subset H_m$ implies $Y^{\le m}(\omega)=Y(\omega)$. Therefore,
for $\omega\in A^-$,
\[
y(\omega)\le t \ \Longrightarrow\ y_U(\omega)\le t
\ \Longleftrightarrow\ Y(\omega)\subset H_t.
\]
Consequently,
\[
\{y\le t\}\subset\{Y\subset H_t\},
\qquad
\PP(y\le t)\le \PP(Y\subset H_t)=C_Y(H_t).
\]

Because $m\le C_Y^{-1}(\alpha)$ and $t<m$, the monotonicity of $C_Y$ implies
$C_Y(H_t)<\alpha$, hence
\[
\PP(y\le t)\ \le\ C_Y(H_t)\ <\ \alpha,\qquad t<m.
\]

\smallskip\noindent
\textbf{Identify $F_y^{-1}(\alpha)$.}
We have shown that $\PP(y\le t)<\alpha$ for all $t<m$, while
$\PP(y\le m)\ge\alpha$. By the definition of the generalized inverse,
\[
F_y^{-1}(\alpha)=\inf\{t:\PP(y\le t)\ge\alpha\}=m.
\]
This completes the proof.
\end{proof}

Taking $\alpha=\tfrac12$, the previous proposition implies that every
point in
\[
[\med_-(y_L),\,\med_+(y_U)]
\]
is a median of some selection of $Y$, so the benchmark median interval
$\Theta_M$ is sharp.

\subsection{Auxiliary rearrangement lemmas}\label{subsec:rearrangement}

We next record two auxiliary results that will be used to characterize
the effect of median restrictions on expectations.

\begin{lemma}[Conditional rearrangement / quantile bound\cite{LiebLoss2001}]\label{lem:rearrangement}
Let $(\Omega,\mathcal F,\PP)$ be a probability space, let $M\in\mathcal F$ with $p_0:=\PP(M)>0$, and let $X:\Omega\to[0,\infty)$ be $\mathcal F$-measurable. Denote by $\mathcal F_M:=\{A\cap M:\ A\in\mathcal F\}$ the trace $\sigma$-algebra on $M$, and by
\[
\PP_M(\,\cdot\,):=\PP(\,\cdot\,\mid M)=\frac{\PP((\,\cdot\,)\cap M)}{p_0}
\]
the conditional probability on $(M,\mathcal F_M)$. Write $F_{X\mid M}(t):=\PP_M(X\le t)$ for the conditional distribution of $X$ given $M$ and let
\[
Q(u):=F^{-1}_{X\mid M}(u):=\inf\{t\in\R:\ F_{X\mid M}(t)\ge u\},\qquad u\in[0,1],
\]
be the (left-continuous) conditional quantile. 

Then,
\textbf{(1)} for any $S\in\mathcal F_M$ with $\PP(S)=s\in[0,p_0]$,
\begin{equation}\label{eq:HL-bound}
\E\big[X\,\mathbf 1_S\big]\ \ge\ p_0\int_0^{s/p_0} Q(u)\,du.
\end{equation} 

\textbf{(2)} If the probability space $(M,\mathcal F_M,\PP_M)$ is atomless, equality in \eqref{eq:HL-bound} is attained by any \emph{least-$X$} set $S$ of mass $s$, i.e.
\[
S=\{X<q\}\cap M\;\cup\;\Big(\{X=q\}\cap M\cap B\Big),
\]
where $q:=Q(s/p_0)$ and $B\in\mathcal F_M$ is chosen with $\PP_M(B)=\theta$ to adjust the measure when $\PP_M(X=q)>0$. 

\textbf{(3)} Without the atomlessness assumption, equality is attained on an atomless product extension by splitting atoms with an auxiliary independent uniform variable.
\end{lemma}

\begin{proof}
\textit{Step 0 (Reduction to the conditional space).}
Since $\PP(S)=s$ if and only if $\PP_M(S)=\alpha:=s/p_0$, and
\[
\E\big[X\,\mathbf 1_S\big] \;=\; p_0\,\E_M\big[X\,\mathbf 1_S\big],
\]
it suffices to prove that, on the probability space $(M,\mathcal F_M,\PP_M)$,
\begin{equation}\label{eq:cond-claim}
\inf_{\substack{S\in\mathcal F_M\\ \PP_M(S)=\alpha}} \E_M[X\,\mathbf 1_S] \;=\; \int_0^{\alpha} Q(u)\,du,
\qquad \alpha\in[0,1].
\end{equation}

\smallskip
\textit{Step 1 (Distributional transform and quantile representation).}
If $(M,\mathcal F_M,\PP_M)$ is not atomless, pass to the product extension
\[
(\bar\Omega,\bar{\mathcal F},\bar{\PP})
\;:=\;
\big(M\times[0,1],\ \mathcal F_M\otimes\mathcal B([0,1]),\ \PP_M\otimes\lambda\big),
\]
where $\lambda$ is Lebesgue measure and let $V:\bar\Omega\to[0,1]$ be the canonical uniform variable. Define on $(\bar\Omega,\bar{\mathcal F},\bar{\PP})$ the \emph{distributional transform}
\[
U\ :=\ F_{X\mid M}(X-)\ +\ V\cdot\Delta F_{X\mid M}(X),
\qquad\text{where}\quad \Delta F_{X\mid M}(x):=F_{X\mid M}(x)-\lim_{t\uparrow x}F_{X\mid M}(t).
\]
Then $U\sim\mathrm{Unif}[0,1]$ (under $\PP_M$ or $\bar\PP$, as appropriate) and
\begin{equation}\label{eq:quantile-rep}
X\ =\ Q(U)\qquad\text{$\PP_M$-a.s.\ (or $\bar\PP$-a.s.)}.
\end{equation}
(When $F_{X\mid M}$ is continuous, one may simply take $U=F_{X\mid M}(X)$ on $M$.)

\smallskip
\textit{Step 2 (From sets to weights depending only on $U$).}
Fix any $S\in\mathcal F_M$ with $\PP_M(S)=\alpha$ and set $H:=\mathbf 1_S$. Using \eqref{eq:quantile-rep} and the tower property,
\[
\E_M[X\,H]\;=\;\E_M\!\big[Q(U)\,H\big]
\;=\;\E_M\!\Big[\E_M\!\big[Q(U)\,H\,\big|\,U\big]\Big]
\;=\;\E_M\!\big[Q(U)\,h(U)\big],
\]
where
\[
h(u)\ :=\ \E_M\!\big[H\,\big|\,U=u\big]\in[0,1]\quad\text{for a.e.\ }u\in[0,1],
\qquad
\int_0^1 h(u)\,du\;=\;\E_M[H]\;=\;\alpha,
\]
because $U\sim\mathrm{Unif}[0,1]$. Thus the minimization in \eqref{eq:cond-claim} reduces to the deterministic calculus of variations problem
\begin{equation}\label{eq:bathtub-program}
\inf\Big\{\ \int_0^1 Q(u)\{h(u)\,du\ :\ 0\le h\le 1,\ \int_0^1 h=\alpha\ \Big\}.
\end{equation}
Here $Q$ is nondecreasing and integrable on $[0,1]$ (since $X\ge0$ is integrable under $\PP_M$).

\smallskip
\textit{Step 3 (Bathtub principle: optimal $h$ is a threshold).}
We show that any minimizer of \eqref{eq:bathtub-program} is of the form
\[
h^*(u)=\mathbf 1_{\{Q(u)<\lambda\}}\;+\;\theta\,\mathbf 1_{\{Q(u)=\lambda\}},
\]
for some $\lambda\in\R$ and $\theta\in[0,1]$ chosen to satisfy $\int_0^1 h^*=\alpha$. The argument is standard:

Suppose $h$ is feasible and not of this form. Then there exist $u_1,u_2$ with $Q(u_1)<Q(u_2)$, $h(u_1)<1$, and $h(u_2)>0$. For $\varepsilon>0$ small, define
\[
\tilde h:=h+\varepsilon\,\mathbf 1_{B_1}-\varepsilon\,\mathbf 1_{B_2},
\]
where $B_1\subset\{u:\ Q(u)\le Q(u_1)\}$ and $B_2\subset\{u:\ Q(u)\ge Q(u_2)\}$ are measurable subsets of equal Lebesgue measure chosen so that $0\le \tilde h\le 1$ a.e.\ and $\int \tilde h=\int h=\alpha$ (possible because $h(u_1)<1$ and $h(u_2)>0$ on sets of positive measure). Then
\[
\int_0^1 Q(u)\,\tilde h(u)\,du
\;=\;\int_0^1 Q(u)\,h(u)\,du
\;+\;\varepsilon\!\left(\int_{B_1} Q(u)\,du - \int_{B_2} Q(u)\,du\right)\]
\[
\;<\;\int_0^1 Q(u)\,h(u)\,du,
\]
since $Q$ is nondecreasing and $Q(u_1)<Q(u_2)$. Hence $h$ is not optimal. Iterating this exchange shows that any optimal $h$ must take values $1$ on the region where $Q$ is smallest, $0$ where $Q$ is largest, and possibly a fractional value on a single level set $\{Q=\lambda\}$ to match the total mass $\alpha$. Because $Q$ is nondecreasing, such an $h^*$ is equivalent (modulo null sets) to the simple threshold
\[
h^*(u)=\mathbf 1_{[0,\alpha]}(u),
\]
when $Q$ is strictly increasing; in general, $h^*$ equals $\mathbf 1_{\{Q<q\}}+\theta\,\mathbf 1_{\{Q=q\}}$ with $q:=Q(\alpha)$ and suitable $\theta\in[0,1]$.

Therefore,
\begin{equation}\label{eq:bathtub-value}
\inf\eqref{eq:bathtub-program}
\;=\;\int_0^1 Q(u)\,h^*(u)\,du
\;=\;\int_0^{\alpha} Q(u)\,du.
\end{equation}

\smallskip
\textit{Step 4 (Attainment by a set $S$).}
If $(M,\mathcal F_M,\PP_M)$ is atomless, there exists $B\in\mathcal F_M$ with any prescribed $\PP_M(B)\in[0,1]$. Let $q:=Q(\alpha)$ and choose $B\subseteq\{X=q\}\cap M$ with $\PP_M(B)=\theta$, where $\theta$ is chosen so that
\[
\PP_M\big(\{X<q\}\cap M\big)+\theta\,\PP_M(\{X=q\}\cap M)=\alpha.
\]
Define
\[
S^* \;:=\; \big(\{X<q\}\cap M\big)\ \cup\ B.
\]
Then $\PP_M(S^*)=\alpha$, and $H^*:=\mathbf 1_{S^*}$ satisfies $H^*=h^*(U)$ in the sense of Step~3; hence
\[
\E_M[X\,\mathbf 1_{S^*}]\;=\;\int_0^{\alpha} Q(u)\,du,
\]
establishing equality in \eqref{eq:cond-claim}.

If $(M,\mathcal F_M,\PP_M)$ has atoms on $\{X=q\}$ and no suitable $B\in\mathcal F_M$ exists with the exact conditional mass, pass to the product extension of Step~1. Use the auxiliary uniform $V$ to split the level set measurably: set
\[
S^* \;:=\; \big(\{X<q\}\cap M\big)\ \cup\ \Big(\{X=q\}\cap M\cap\{V\le \theta\}\Big),
\]
now a set in $\bar{\mathcal F}$ with $\bar{\PP}(S^*)=\alpha$ and $\E_{\bar\PP}[X\,\mathbf 1_{S^*}]=\int_0^{\alpha}Q(u)\,du$ (note $X$ does not depend on $V$). This yields equality on the extension.

\smallskip
\textit{Step 5.}
Combining Steps~1–4, we have proved \eqref{eq:cond-claim}. Multiplying by $p_0$ gives \eqref{eq:HL-bound}. The stated forms of the equality sets follow from the constructions above.

\smallskip
\textit{Edge cases.}
If $\alpha=0$ (i.e., $s=0$), both sides vanish. If $\alpha=1$ (i.e., $s=p_0$), both sides equal $\E[X\,\mathbf 1_M]=p_0\int_0^1 Q(u)\,du$.
\end{proof}

\begin{lemma}[Quantile--area identity \cite{HardyLittlewoodPolya1952}]\label{lem:quantile_area}
Let $Z$ be a nonnegative real-valued random variable with distribution function $F_Z$ and
(left-continuous) quantile function
\[
F_Z^{-1}(u) \;:=\; \inf\{t\in\mathbb{R}: F_Z(t)\ge u\},\qquad u\in(0,1).
\]
Then for every $\alpha\in(0,1)$,
\begin{equation}\label{eq:quantile-area-lemma}
  \int_0^\alpha F_Z^{-1}(u)\,du
  \;=\;
  \int_0^{\infty} (\alpha - F_Z(t))_+\,dt,
\end{equation}
where $(x)_+ := \max\{x,0\}$.
\end{lemma}

\begin{proof}
Set $Q(u):=F_Z^{-1}(u)$ for $u\in(0,1)$ and fix $\alpha\in(0,1)$. Since $Z\ge 0$ almost surely, we have $Q(u)\ge 0$ for all $u\in(0,1)$. By the layer-cake representation and Tonelli's theorem,
\begin{align*}
\int_0^\alpha Q(u)\,du
&= \int_0^\alpha \left(\int_0^\infty \mathbf{1}_{\{t<Q(u)\}}\,dt\right)du \\
&= \int_0^\infty \left(\int_0^\alpha \mathbf{1}_{\{t<Q(u)\}}\,du\right)dt.
\end{align*}
For the left-continuous quantile function
\[
Q(u)=\inf\{x\in\mathbb R:F_Z(x)\ge u\},
\]
we have, for every $t\ge0$ and $u\in(0,1)$,
\[
t<Q(u)\quad\Longleftrightarrow\quad F_Z(t)<u.
\]
Hence, for fixed $t\ge0$,
\[
\int_0^\alpha \mathbf{1}_{\{t<Q(u)\}}\,du
=
\int_0^\alpha \mathbf{1}_{\{F_Z(t)<u\}}\,du
=
(\alpha-F_Z(t))_+.
\]
Substituting this into the previous display gives
\[
\int_0^\alpha Q(u)\,du
=
\int_0^\infty (\alpha-F_Z(t))_+\,dt,
\]
which is exactly \eqref{eq:quantile-area-lemma}.
\end{proof}
\subsection{Bounds on \texorpdfstring{$\E(y)$}{E(y)} over \texorpdfstring{$\Sel(Y\mid m)$}{Sel(Y|m)}}\label{subsec:IdMean}

We now consider the effect of knowing a median on the range of admissible
expectations. Let $Y=[y_L,y_U]$ be a random interval with
$\E|y_L|+\E|y_U|<\infty$, and consider the restricted selection set $\Sel(Y\mid m)$ of all selections with median $m$ in the sense that
\[
\PP(y\le m)\ge\tfrac12,\qquad \PP(y\ge m)\ge\tfrac12.
\]
We are interested in the range
\[
\big\{\E[y]: y\in\Sel(Y\mid m)\big\}.
\]

Define
\[
A_-:=\{y_U<m\},\qquad A_+:=\{y_L>m\},\qquad M:=\{y_L\le m\le y_U\},
\]
with $p_-:=\PP(A_-)$, $p_+:=\PP(A_+)$, $p_0:=\PP(M)=1-p_--p_+$. Set the required “mass shortfalls’’
\[
\alpha_-(m):=\big(\tfrac12-p_-\big)_+,\qquad \alpha_+(m):=\big(\tfrac12-p_+\big)_+,
\]
and on $M$ define the nonnegative gaps to the median
\[
U_m:=y_U-m,\qquad L_m:=m-y_L.
\]
Let $F^{-1}_{U_m\mid M}$ and $F^{-1}_{L_m\mid M}$ denote the (left-continuous) quantile functions of the conditional laws of $U_m$ and $L_m$ given $M$.

\begin{proposition}[Median-restricted Aumann expectation]\label{prop:MedianRestrictedAumann}
Assume \emph{feasibility}, $p_-\le\tfrac12$ and $p_+\le\tfrac12$. The attainable set of means under the median restriction
\[
\Sel(Y\mid m)=\{y\in\Sel(Y): m\in\Med(y)\}.
\]
is the closed interval
\begin{equation}\label{eq:MedianMeanID}
\{\E[y]:y\in\mathcal S_m\}\;=\;\Big[\,\E_{\min}(m),\ \E_{\max}(m)\,\Big],
\end{equation}
with endpoints
\begin{equation}\label{eq:MedianMeanEndpoints}
\E_{\max}(m)\;=\;\E[y_U]\;-\;p_0\int_{0}^{\alpha_-(m)/p_0}\! F^{-1}_{\,U_m\mid M}(u)\,du,
\end{equation}
and
\begin{equation}
\E_{\min}(m)\;=\;\E[y_L]\;+\;p_0\int_{0}^{\alpha_+(m)/p_0}\! F^{-1}_{\,L_m\mid M}(u)\,du,
\end{equation}
where the integrals are understood as $0$ when $p_0=0$. Both endpoints are attained by measurable selections, and every value in \([\E_{\min}(m),\E_{\max}(m)]\) is attainable.
\end{proposition}

\begin{proof}
Recall
\[
A_-:=\{y_U<m\},\quad A_+:=\{y_L>m\},\quad M:=\{y_L\le m\le y_U\},
\]
with $p_-:=\PP(A_-)$, $p_+:=\PP(A_+)$, $p_0:=\PP(M)$, the shortfalls
\[
\alpha_-(m):=(\tfrac12-p_-)_+,\qquad \alpha_+(m):=(\tfrac12-p_+)_+,
\]
and the nonnegative gaps on $M$,
\[
U_m:=y_U-m,\qquad L_m:=m-y_L.
\]
We assume feasibility $p_-\le\tfrac12$ and $p_+\le\tfrac12$, which implies $\alpha_-(m)+\alpha_+(m)=p_0$ and ensures that the median constraint can be satisfied by allocating on $M$ a mass $\alpha_-(m)$ at or below $m$ and a mass $\alpha_+(m)$ at or above $m$.

\medskip\noindent
\emph{Upper endpoint.}
Let $y\in\Sel(Y)$ satisfy the median restriction $\PP(y\le m)\ge\tfrac12$. Since $y\le m$ on $A_-$ and $y\le y_U$ on $M$, we must have
\[
\PP\big(\{y\le m\}\cap M\big)\ \ge\ \tfrac12-\PP(A_-)\ =\ \alpha_-(m).
\]
On $\{y\le m\}\cap M$, $y_U-y\ge U_m$. Therefore
\[
\E[y_U-y]\ \ge\ \E\!\big[U_m\,\1_{\{y\le m\}\cap M}\big]\ \ge\ \inf_{S\subset M:\ \PP(S)\ge\alpha_-(m)} \E\!\big[U_m\,\1_{S}\big].
\]
By Lemma~\ref{lem:rearrangement} applied to $X=U_m\ge0$ on $M$,
\[
\inf_{S\subset M:\ \PP(S)\ge\alpha_-(m)} \E\!\big[U_m\,\1_{S}\big]\ =\ p_0\int_0^{\alpha_-(m)/p_0} F^{-1}_{U_m\mid M}(u)\,du.
\]
Consequently,
\[
\E[y]\ \le\ \E[y_U]-\E[y_U-y]\ \le\ \E[y_U]\;-\;p_0\int_0^{\alpha_-(m)/p_0} F^{-1}_{U_m\mid M}(u)\,du\ :=\ \E_{\max}(m).
\]

\medskip\noindent
\emph{Lower endpoint.}
Symmetrically, the median restriction implies
\[
\PP\big(\{y\ge m\}\cap M\big)\ \ge\ \tfrac12-\PP(A_+)\ =\ \alpha_+(m),
\]
and on $\{y\ge m\}\cap M$ one has $y-y_L\ge L_m$. Thus
\[
\E[y-y_L]\ \ge\ \inf_{S\subset M:\ \PP(S)\ge\alpha_+(m)} \E\!\big[L_m\,\1_{S}\big]
\ =\ p_0\int_0^{\alpha_+(m)/p_0} F^{-1}_{L_m\mid M}(u)\,du,
\]
by Lemma~\ref{lem:rearrangement}, which yields
\[
\E[y]\ \ge\ \E[y_L]\;+\;p_0\int_0^{\alpha_+(m)/p_0} F^{-1}_{L_m\mid M}(u)\,du\ :=\ \E_{\min}(m).
\]

\medskip\noindent
\emph{Attainment.}
Let $q_-:=F^{-1}_{U_m\mid M}(\alpha_-(m)/p_0)$ and choose a measurable $S_-\subset M$ with $\PP(S_-)=\alpha_-(m)$ such that $U_m$ is $\PP$-a.s.\ minimal on $S_-$, i.e.\ $S_-=\{U_m<q_-\}\cup(\{U_m=q_-\}\cap B_-)$ with suitable measurable $B_-$. Define
\[
\bar y(\omega):=\begin{cases}
m, & \omega\in S_-,\\
y_U(\omega), & \omega\in M\setminus S_-,\\
y_U(\omega), & \omega\in A_-\cup A_+.
\end{cases}
\]
Then $\bar y\in\Sel(Y)$.  We verify the median condition $m\in\Med(\bar y)$.

On $A_-$ we have $\bar y = y_U < m$ (since $A_- = \{y_U < m\}$); on $S_- \subset M$ we have $\bar y = m$; and on $M \setminus S_-$ we have $\bar y = y_U \ge m$ (with equality when $y_U = m$).  On $A_+$ we have $\bar y = y_U \ge y_L > m$.  Therefore
\[
\PP(\bar y \le m)
\;=\;
\PP(A_-) + \PP(S_-) + \PP\bigl(\{y_U = m\} \cap (M \setminus S_-)\bigr)
\;\ge\;
p_- + \alpha_-(m)
\;=\;
\tfrac12.
\]
For the other half: on $A_-$ we have $\bar y = y_U < m$, so $\bar y \ge m$ fails there; everywhere else $\bar y \ge m$.  Hence
\[
\PP(\bar y \ge m)
\;=\;
1 - \PP(A_-)
\;=\;
1 - p_-
\;\ge\;
\tfrac12,
\]
using the feasibility assumption $p_- \le \tfrac12$.  Thus $m \in \Med(\bar y)$, and by construction
\[
\E[y_U-\bar y]\ =\ \E\!\big[U_m\,\1_{S_-}\big]\ =\ p_0\int_0^{\alpha_-(m)/p_0} F^{-1}_{U_m\mid M}(u)\,du,
\]
so $\E[\bar y]=\E_{\max}(m)$.

Likewise, let $q_+:=F^{-1}_{L_m\mid M}(\alpha_+(m)/p_0)$ and choose a measurable $S_+\subset M$ of measure $\alpha_+(m)$ which minimizes $L_m$. Define
\[
\underline y(\omega):=\begin{cases}
m & \omega\in S_+,\\
y_L(\omega), & \omega\in M\setminus S_+,\\
y_L(\omega), & \omega\in A_-\cup A_+.
\end{cases}
\]
Then $\underline y\in\Sel(Y)$.  On $A_+$ we have $\underline y = y_L > m$; on $S_+ \subset M$ we have $\underline y = m$; and on $M \setminus S_+$ we have $\underline y = y_L \le m$ (with equality when $y_L = m$).  On $A_-$ we have $\underline y = y_L \le y_U < m$.  Therefore
\[
\PP(\underline y \ge m)
\;=\;
\PP(A_+) + \PP(S_+) + \PP\bigl(\{y_L = m\} \cap (M \setminus S_+)\bigr)
\;\ge\;
p_+ + \alpha_+(m)
\;=\;
\tfrac12,
\]
and
\[
\PP(\underline y \le m)
\;=\;
1 - \PP(A_+)
\;=\;
1 - p_+
\;\ge\;
\tfrac12.
\]
Thus $m \in \Med(\underline y)$, and
\[
\E[\underline y-y_L]\ =\ \E\!\big[L_m\,\1_{S_+}\big]\ =\ p_0\int_0^{\alpha_+(m)/p_0} F^{-1}_{L_m\mid M}(u)\,du,
\]
so $\E[\underline y]=\E_{\min}(m)$.

\medskip\noindent
\emph{Filled interval.}
For $\theta\in[0,1]$ put $y_\theta:=\theta\,\bar y+(1-\theta)\,\underline y$. Since $Y(\omega)$ is an interval, $y_\theta(\omega)\in Y(\omega)$ a.s., hence $y_\theta\in\Sel(Y)$. Moreover,
\[
\{y_\theta\le m\}\supseteq A_-\cup S_-,\qquad \{y_\theta\ge m\}\supseteq A_+\cup S_+,
\]
so $\PP(y_\theta\le m)\ge\tfrac12$ and $\PP(y_\theta\ge m)\ge\tfrac12$; thus $y_\theta\in\mathcal S_m$. Linearity of expectation gives
\[
\E[y_\theta]=\theta\,\E_{\max}(m)+(1-\theta)\,\E_{\min}(m),
\]
so every value in $[\E_{\min}(m),\E_{\max}(m)]$ is attained. This establishes \eqref{eq:MedianMeanID}–\eqref{eq:MedianMeanEndpoints}.
\end{proof}

\begin{remark}[Remarks on edge cases]
If $p_0=0$ and feasibility holds, one must have $p_-=p_+=\tfrac12$, the integrals in \eqref{eq:MedianMeanEndpoints} vanish, and the median restriction does not tighten $\Aum[Y]=[\E y_L,\E y_U]$. If $m\in [y_L,y_U]$ a.s., then $p_-=p_+=0$, $\alpha_-(m)=\alpha_+(m)=\tfrac12$, and both quantile integrals run to $1/2$, recovering the familiar constant-threshold formulas (e.g., for $Y=[a,b]$ a.s., the interval $[\tfrac{a+m}{2},\tfrac{m+b}{2}]$).
\end{remark}

\begin{figure}[t]
    \centering
    \includegraphics[width=0.65\textwidth]{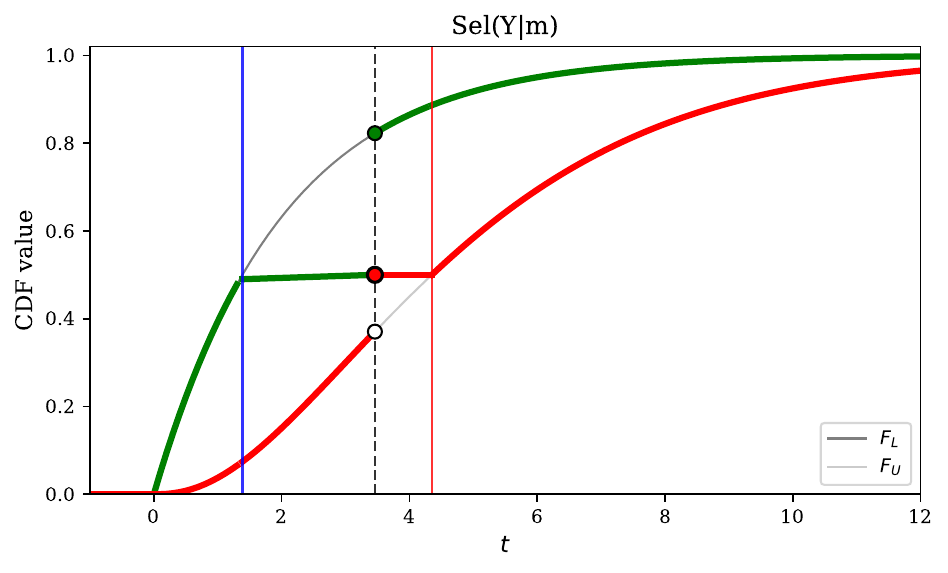}
      \caption{\small{Schematic construction of the extremal boundary distributions in Example~\ref{ex:schematic}. Thin black and gray curves show the
    baseline CDFs $F_L$ and $F_U$ of the chi-square lower and upper bounds.  The bold green
    (resp.\ red) segments indicate the local modification of $F_L$ on $[m_L,m]$
    (resp.\ of $F_U$ on $[m,m_U]$) required to satisfy the median constraints
    $\PP(y\le m)=\PP(y\ge m)=\tfrac12$ for an underlying selection $y\in\Sel(Y)$.
    Vertical lines at $m_L$, $m$, and $m_U$ and the bullets at height $1/2$
    highlight the binding median restrictions.}}
    \label{fig:median-mean-boundary}
\end{figure}

\begin{example}[Median-restricted Aumann expectation]\label{ex:schematic}
To illustrate Proposition~\ref{prop:MedianRestrictedAumann}, consider a stylized selection model
\[
Y=[y_L,y_U],
\]
where $y_L$ and $y_U$ are constructed as follows. Let $U\sim\mathrm{Unif}(0,1)$
and write $F_L,F_U$ for the distribution functions of $\chi^2_{(2)}$ and $\chi^2_{(5)}$.
Set
\[
y_L := F_L^{-1}(U), \qquad y_U := F_U^{-1}(U).
\]
Then $y_L\sim\chi^2_{(2)}$, $y_U\sim\chi^2_{(5)}$, and $y_L\le y_U$ almost surely.

The medians of the lower and upper bounds are
\[
m_L := \Med(y_L) = 1.386,\qquad
m_U := \Med(y_U) = 4.351,
\]
and we fix a target median
\[
m := 0.3\,m_L + 0.7\,m_U = 3.46.
\]
We fix $m$ as the median constraint, so that $\Sel(Y\mid m)$ consists of all $y\in\Sel(Y)$ satisfying that $\PP(y\le m)\ge\tfrac12$ and $\PP(y\ge m)\ge\tfrac12$.

In this configuration the bounded region for the mean simplifies to
\[
\E[y] \in \big[\,\E(y_L)+S_L,\ \E(y_U)-S_U\,\big] \qquad\text{for all }y\in\Sel(Y\mid m),
\]
where
\[
S_L = \int_{m_L}^{m}\!\big(F_L(t)-\tfrac12\big)\,dt,\qquad
S_U = \int_{m}^{m_U}\!\big(\tfrac12-F_U(t)\big)\,dt.
\]
The extremal CDFs that achieve the endpoints $\E_{\min}(m)$ and $\E_{\max}(m)$ in
\eqref{eq:MedianMeanID}–\eqref{eq:MedianMeanEndpoints} are obtained by locally reshaping $F_L$
on $[m_L,m]$ and $F_U$ on $[m,m_U]$ so that exactly one half of the mass of the induced selection
lies below and above $m$, respectively.  This construction is visualized in
Figure~\ref{fig:median-mean-boundary}.
\end{example}

\begin{remark}[Derivation of the bounds in
  Example~\ref{ex:schematic}]\label{rem:SL_SU_derivation}
We show how to pass from the general formulas of
Proposition~\ref{prop:MedianRestrictedAumann} to explicit CDF
integrals.  The key tool is the quantile--area identity
(Lemma~\ref{lem:quantile_area}): for any nonnegative random variable $Z$
with CDF $F_Z$ and any $\alpha\in(0,1)$,
\begin{equation}\label{eq:quantile-area}
  \int_0^\alpha F_Z^{-1}(u)\,du
  \;=\;
  \int_{0}^{\infty} (\alpha - F_Z(t))_+\,dt.
\end{equation}

Write $p_-=\PP(y_U<m)$, $p_+=\PP(y_L>m)$,
$p_0=\PP(y_L\le m\le y_U)=1-p_--p_+$, and
$\alpha_-=\tfrac12-p_-$, $\alpha_+=\tfrac12-p_+$ (assuming
feasibility $p_-,p_+\le\tfrac12$).

\medskip\noindent\textbf{Upper endpoint.}\quad
By Proposition~\ref{prop:MedianRestrictedAumann},
\[
\E_{\max}(m)
\;=\;
\E[y_U]
\;-\;
p_0\!\int_0^{\,\alpha_-/p_0}\!F^{-1}_{U_m\mid M}(u)\,du,
\]
where $U_m=y_U-m$ and the quantile is conditional on
$M=\{y_L\le m\le y_U\}$.  Applying identity
\eqref{eq:quantile-area} on the conditional space
$(M,\PP_M)$ with $Z=U_m$ and $\alpha=\alpha_-/p_0$ gives
\[
p_0\!\int_0^{\alpha_-/p_0}\!F^{-1}_{U_m\mid M}(u)\,du
\;=\;
p_0\!\int_0^{\infty}
\Bigl(\frac{\alpha_-}{p_0}-F_{U_m\mid M}(\ell)\Bigr)_+d\ell.
\]
Since $F_{U_m\mid M}(\ell)
=\PP_M(y_U-m\le\ell)
=\PP_M(y_U\le m+\ell)
=\frac{\PP(y_U\le m+\ell,\;y_L\le m\le y_U)}{p_0}$,
and writing $G(\ell):=\PP(y_U\le m+\ell,\;y_L\le m\le y_U)$ so
that $F_{U_m\mid M}(\ell)=G(\ell)/p_0$, we obtain
\[
p_0\!\int_0^{\alpha_-/p_0}\!F^{-1}_{U_m\mid M}(u)\,du
\;=\;
\int_0^{\infty}
\bigl(\alpha_- - G(\ell)\bigr)_+\,d\ell.
\]
In the special case $p_-=p_+=0$ (equivalently $m\in[y_L,y_U]$ a.s.),
we have $p_0=1$, $\alpha_-=\alpha_+=\tfrac12$, and
$G(\ell)=F_U(m+\ell)$, so the integral reduces to
\[
\int_0^{\infty}\bigl(\tfrac12-F_U(m+\ell)\bigr)_+\,d\ell
\;=\;
\int_m^{m_U}\bigl(\tfrac12-F_U(t)\bigr)\,dt
\;=:\;S_U.
\]
The analogous calculation for the lower endpoint yields
$S_L=\int_{m_L}^{m}(F_L(t)-\tfrac12)\,dt$ under the same
simplifying assumption.

\medskip
In the chi-squared configuration of Example~\ref{ex:schematic},
$p_-$ and $p_+$ are strictly positive, so the simple formulas
involving $S_L$ and $S_U$ are only approximate.  However, the
qualitative picture --- the shaded areas between the CDFs and the
$1/2$--level being ``costs'' that shrink the Aumann range ---
remains valid.  For exact numerical evaluation one should use
the general conditional formulas of
Proposition~\ref{prop:MedianRestrictedAumann} directly.
\end{remark}

\subsection{Attainable medians under a mean restriction}

We now consider the reverse question: does the mean constraint $\E[y]=\kappa$ necessarily sharpen the benchmark median interval
$\Theta_M=[\Med(y_L),\Med(y_U)]$? The next example shows that the answer is
negative in general.

\begin{example}[Mean restriction need not shrink the attainable median set]\label{ex:no-shrink-med}
Fix $\kappa\in\R$ and $d>0$. Let $(\Omega,\mathcal A,\PP)$ be a non-atomic probability space, and let
$B\in\mathcal A$ satisfy $\PP(B)=\tfrac12$. Define the random interval
\[
Y(\omega)=
\begin{cases}
[\kappa-2d,\ \kappa], & \omega\in B,\\[2mm]
[\kappa,\ \kappa+2d], & \omega\in B^c.
\end{cases}
\]
Equivalently,
\[
(y_L,y_U)=
\begin{cases}
(\kappa-2d,\ \kappa), & \omega\in B,\\
(\kappa,\ \kappa+2d), & \omega\in B^c.
\end{cases}
\]

We first identify the benchmark median interval. Since
\[
y_L=
\begin{cases}
\kappa-2d, & \omega\in B,\\
\kappa, & \omega\in B^c,
\end{cases}
\qquad
y_U=
\begin{cases}
\kappa, & \omega\in B,\\
\kappa+2d, & \omega\in B^c,
\end{cases}
\]
each with probabilities $1/2$ and $1/2$, we have
\[
\Med(y_L)=[\kappa-2d,\kappa],
\qquad
\Med(y_U)=[\kappa,\kappa+2d].
\]
Hence the benchmark interval from \eqref{equ:benchmarkBounds_M} is
\[
\Theta_M=[\kappa-2d,\kappa+2d].
\]

Now fix any target median
\[
m\in[\kappa-2d,\kappa+2d].
\]
We show that there exists a selection $y\in\Sel(Y)$ such that
\[
\E[y]=\kappa
\qquad\text{and}\qquad
m\in\Med(y).
\]

\smallskip
\noindent
\emph{Case 1: $m\in[\kappa-2d,\kappa]$.}
Define
\[
y(\omega)=
\begin{cases}
m, & \omega\in B,\\
2\kappa-m, & \omega\in B^c.
\end{cases}
\]
Since
\[
m\in[\kappa-2d,\kappa]
\qquad\text{and}\qquad
2\kappa-m\in[\kappa,\kappa+2d],
\]
we have $y(\omega)\in Y(\omega)$ almost surely, so $y\in\Sel(Y)$. Moreover,
\[
\E[y]=\tfrac12\,m+\tfrac12\,(2\kappa-m)=\kappa.
\]
Because $m\le\kappa\le 2\kappa-m$, the distribution of $y$ puts probability
$1/2$ at $m$ and probability $1/2$ at $2\kappa-m$, hence
\[
\PP(y\le m)=\tfrac12,
\qquad
\PP(y\ge m)=1.
\]
Therefore $m\in\Med(y)$.

\smallskip
\noindent
\emph{Case 2: $m\in[\kappa,\kappa+2d]$.}
Define
\[
y(\omega)=
\begin{cases}
2\kappa-m, & \omega\in B,\\
m, & \omega\in B^c.
\end{cases}
\]
Now
\[
2\kappa-m\in[\kappa-2d,\kappa]
\qquad\text{and}\qquad
m\in[\kappa,\kappa+2d],
\]
so again $y\in\Sel(Y)$. Also,
\[
\E[y]=\tfrac12(2\kappa-m)+\tfrac12 m=\kappa.
\]
Since $2\kappa-m\le \kappa\le m$, we get
\[
\PP(y\le m)=1,
\qquad
\PP(y\ge m)=\tfrac12,
\]
and therefore $m\in\Med(y)$.

\smallskip
Since $m$ was arbitrary in $[\kappa-2d,\kappa+2d]$, we conclude that
\[
\{\,m\in\R:\exists\,y\in\Sel(Y)\text{ with }\E[y]=\kappa\text{ and }m\in\Med(y)\,\}
=
[\kappa-2d,\kappa+2d]
=
\Theta_M.
\]
Thus, in this random-interval design, the mean restriction $\E[y]=\kappa$
does \emph{not} shrink the benchmark attainable median interval.
\end{example}
In particular, the mean restriction need not reduce the \emph{set} of attainable medians,
even though it may reduce the set of attainable distributions.
\subsection{Event probabilities: unrestricted and mean-restricted ranges}

We now turn to the range of probabilities $\PP(y\in A)$ as $y$ varies over selections of $Y$, with and without a mean restriction. Let $A$ be a fixed
Borel subset of $\R$.

\begin{proposition}[Unrestricted random-set bounds]\label{prop:RS-unrestricted}
Let $A$ be a closed subset of $\R$. Then for any $y\in\Sel(Y)$,
\begin{equation}\label{eq:RS-bounds}
\PP\!\big(Y\subseteq A\big)\ \le\ \PP(y\in A)\ \le\ \PP\!\big(Y\cap A\neq\varnothing\big).
\end{equation}
Both bounds are attained by measurable selections.
\end{proposition}

\begin{proof}
If $Y(\omega)\subseteq A$ and $y(\omega)\in Y(\omega)$, then $y(\omega)\in A$; thus
$\{Y\subseteq A\}\subseteq\{y\in A\}$ and
$\PP(y\in A)\ge\PP(Y\subseteq A)$.
Similarly, if $y(\omega)\in A$ then $Y(\omega)\cap A\neq\varnothing$, so
$\{y\in A\}\subseteq\{Y\cap A\neq\varnothing\}$ and $\PP(y\in A)\le\PP(Y\cap A\neq\varnothing)$.

For sharpness, define
\[
A_1:=\{\omega:Y(\omega)\subseteq A\},\qquad
A_2:=\{\omega:Y(\omega)\cap A\neq\varnothing\}.
\]

For the upper bound, define a set-valued map
\[
G_+(\omega):=
\begin{cases}
Y(\omega)\cap A, & \omega\in A_2,\\
Y(\omega), & \omega\notin A_2.
\end{cases}
\]
This is a measurable random closed set with nonempty values, so by measurable
selection there exists $y^+$ with $y^+(\omega)\in G_+(\omega)$ almost surely
(see \cite{CastaingValadier1977,Molchanov2005}). Then $y^+(\omega)\in A$
whenever $Y(\omega)\cap A\neq\emptyset$, and necessarily $y^+(\omega)\notin A$
when $Y(\omega)\cap A=\emptyset$. Hence
\[
\PP(y^+\in A)=\PP(Y\cap A\neq\emptyset),
\]
so the upper bound in \eqref{eq:RS-bounds} is attainable.

For the lower bound, define
\[
G_-(\omega):=
\begin{cases}
Y(\omega)\setminus A, & \omega\notin A_1,\ \text{if this set is nonempty},\\
Y(\omega), & \omega\in A_1.
\end{cases}
\]
Again $G_-$ is a measurable random closed set with nonempty values. A measurable
selection $y^-$ with $y^-(\omega)\in G_-(\omega)$ almost surely exists. By
construction $y^-(\omega)\notin A$ whenever $Y(\omega)\not\subseteq A$, and
$y^-(\omega)\in A$ if $Y(\omega)\subseteq A$. Thus
\[
\PP(y^-\in A)=\PP(Y\subseteq A),
\]
and the lower bound is attained.
\end{proof}

We next impose the mean restriction and define, for $\kappa\in[\EE y_L,\EE y_U]$,
\[
\Sel(Y\,|\,\kappa):=\{\,y\in\Sel(Y):\ \EE[y]=\kappa\,\},
\]
and
\[
U(\kappa):=\sup_{y\in\Sel(Y\,|\,\kappa)}\PP(y\in A),\quad
L(\kappa):=\inf_{y\in\Sel(Y\,|\,\kappa)}\PP(y\in A).
\]
Introduce the $\omega$-wise extremal points of $A$ inside $Y$ and the associated gaps
\begin{align}
a^+(\omega)&:=\sup\big(Y(\omega)\cap A\big)\in[-\infty,\infty],&
a^-(\omega)&:=\inf\big(Y(\omega)\cap A\big)\in[-\infty,\infty], \label{eq:a+-}\\
\Delta^+(\omega)&:=\big(y_U(\omega)-a^+(\omega)\big)_+,&
\Delta^-(\omega)&:=\big(a^-(\omega)-y_L(\omega)\big)_+,
\end{align}
with the conventions $\sup\varnothing=-\infty$ and $\inf\varnothing=+\infty$.

The next result uses a Lagrangian and an interchange of $\sup$ and expectation
in the spirit of \cite{Artstein1974,Molchanov2005}; convex duality arguments are standard, see e.g.\ \cite{Rockafellar1970}.

\begin{theorem}[Dual characterization under a mean restriction]\label{thm:dual}
Define, for $\lambda\in\R$ and $\omega\in\Omega$,
\[
\Psi(\lambda,\omega)=
\begin{cases}
\max\{\lambda y_U(\omega),\ 1+\lambda a^+(\omega)\}, & \lambda\ge 0,\\[2pt]
\max\{\lambda y_L(\omega),\ 1+\lambda a^-(\omega)\}, & \lambda< 0,
\end{cases}
\]
and
\[
\Phi(\lambda,\omega)=
\begin{cases}
\min\{\lambda y_L(\omega),\ 1+\lambda a^-(\omega)\}, & \lambda\ge 0,\\[2pt]
\min\{\lambda y_U(\omega),\ 1+\lambda a^+(\omega)\}, & \lambda< 0.
\end{cases}
\]
Then
\begin{equation}\label{eq:dual-envelopes-kappa}
U(\kappa)=\inf_{\lambda\in\R}\Big\{\EE\big[\Psi(\lambda,\cdot)\big]-\lambda\kappa\Big\},\qquad
L(\kappa)=\sup_{\lambda\in\R}\Big\{\EE\big[\Phi(\lambda,\cdot)\big]-\lambda\kappa\Big\}.
\end{equation}
Moreover, $U$ is concave in $\kappa$, while $L$ is convex. 
\end{theorem}

\begin{proof}
We prove the formula for $U(\kappa)$; the argument for $L(\kappa)$ is analogous.

For each $\lambda\in\R$, consider the Lagrangian
\[
\mathcal L(y,\lambda)
:=
\E\big[\mathbf 1_A(y)\big]+\lambda\big(\E[y]-\kappa\big)
=
\E\big[\mathbf 1_A(y)+\lambda y\big]-\lambda\kappa,
\qquad y\in\Sel(Y).
\]
For fixed $\lambda$, maximizing $\mathcal L(y,\lambda)$ over $y\in\Sel(Y)$ is
equivalent to maximizing
\[
\E\big[f_\lambda(\omega,y(\omega))\big],
\qquad
f_\lambda(\omega,x):=\mathbf 1_A(x)+\lambda x,
\]
subject to the pointwise constraint $y(\omega)\in Y(\omega)$. By the standard
interchange theorem for random closed sets (see, e.g.,
\cite{Artstein1974,Molchanov2005}),
\[
\sup_{y\in\Sel(Y)} \E[f_\lambda(\omega,y(\omega))]
=
\E\Big[\sup_{x\in Y(\omega)} f_\lambda(\omega,x)\Big].
\]

If $\lambda\ge 0$, then for each $\omega$,
\[
\sup_{x\in Y(\omega)} f_\lambda(\omega,x)
=
\max\Big\{
\sup_{x\in Y(\omega)\setminus A}\lambda x,\,
\sup_{x\in Y(\omega)\cap A}(1+\lambda x)
\Big\}
=
\max\{\lambda y_U(\omega),\,1+\lambda a^+(\omega)\},
\]
with the convention that the second term is $-\infty$ when
$Y(\omega)\cap A=\varnothing$. If $\lambda<0$, similarly,
\[
\sup_{x\in Y(\omega)} f_\lambda(\omega,x)
=
\max\{\lambda y_L(\omega),\,1+\lambda a^-(\omega)\}.
\]
Hence
\[
\sup_{y\in\Sel(Y)}\mathcal L(y,\lambda)
=
\E[\Psi(\lambda,\cdot)]-\lambda\kappa.
\]

On the other hand, if $y\in\Sel(Y\mid\kappa)$, then $\E[y]=\kappa$, so
\[
\mathcal L(y,\lambda)=\PP(y\in A)
\qquad\text{for every }\lambda\in\R.
\]
Therefore
\[
U(\kappa)
=
\sup_{y\in\Sel(Y\mid\kappa)}\PP(y\in A)
=
\sup_{y\in\Sel(Y)}\inf_{\lambda\in\R}\mathcal L(y,\lambda).
\]

It remains to justify the minimax interchange. Define the feasible set
\[
\mathcal Y
:=
\bigl\{(\E[y],\,\E[\mathbf 1_A(y)]):\,y\in\Sel(Y)\bigr\}
\subset\R^2.
\]
Then
\[
U(\kappa)
=
\sup\{\,p:(\kappa,p)\in\mathcal Y\,\}
=
\sup_{(m,p)\in\mathcal Y}\inf_{\lambda\in\R}\bigl[p+\lambda(m-\kappa)\bigr].
\]

We first show that $\mathcal Y$ is convex. Let
\[
(m_i,p_i)\in\mathcal Y,\qquad i=1,2,
\]
so there exist selections $y_i\in\Sel(Y)$ with
\[
m_i=\E[y_i],\qquad p_i=\E[\mathbf 1_A(y_i)].
\]
Fix $\theta\in[0,1]$. Pass to the atomless product extension
\[
(\bar\Omega,\bar{\mathcal A},\bar\PP)
:=
(\Omega\times[0,1],\,\mathcal A\otimes\mathcal B([0,1]),\,\PP\otimes\lambda),
\]
where $\lambda$ is Lebesgue measure, and let $U$ be the coordinate uniform
random variable on $[0,1]$, independent of all random objects on $\Omega$.
Extend $Y,y_1,y_2$ to $\bar\Omega$ by ignoring the second coordinate, and define
\[
\bar y(\omega,u)
:=
\mathbf 1_{\{u\le\theta\}}\,y_1(\omega)
+
\mathbf 1_{\{u>\theta\}}\,y_2(\omega).
\]
Then $\bar y(\omega,u)\in Y(\omega)$ for all $(\omega,u)$, so
$\bar y\in\Sel(Y)$ on the product space. By independence of $U$,
\[
\E_{\bar\PP}[\bar y]
=
\theta\,\E[y_1]+(1-\theta)\,\E[y_2]
=
\theta m_1+(1-\theta)m_2,
\]
and
\[
\E_{\bar\PP}[\mathbf 1_A(\bar y)]
=
\theta\,\E[\mathbf 1_A(y_1)]
+
(1-\theta)\,\E[\mathbf 1_A(y_2)]
=
\theta p_1+(1-\theta)p_2.
\]
Thus
\[
\theta(m_1,p_1)+(1-\theta)(m_2,p_2)\in\mathcal Y,
\]
after identifying the product extension with the original atomless setting.
Hence $\mathcal Y$ is convex.

Next, $\mathcal Y$ is bounded, because every selection satisfies
$y_L\le y\le y_U$ a.s., so
\[
\E[y]\in[\E y_L,\E y_U],
\qquad
\E[\mathbf 1_A(y)]\in[0,1].
\]
Thus
\[
\mathcal Y\subset[\E y_L,\E y_U]\times[0,1].
\]

We now note that the upper boundary function
\[
\phi(m):=\sup\{\,p:(m,p)\in\mathcal Y\,\},
\qquad
m\in[\E y_L,\E y_U],
\]
is concave, since $\mathcal Y$ is convex. In particular,
\[
U(\kappa)=\phi(\kappa).
\]
A basic fact from finite-dimensional convex analysis is that every finite
concave function on an interval is the pointwise infimum of its affine majorants.
Equivalently,
\[
\phi(\kappa)
=
\inf_{\lambda\in\R}\sup_{(m,p)\in\mathcal Y}\bigl[p+\lambda(m-\kappa)\bigr].
\]
Therefore
\[
U(\kappa)
=
\inf_{\lambda\in\R}\sup_{(m,p)\in\mathcal Y}\bigl[p+\lambda(m-\kappa)\bigr].
\]
Since each $(m,p)\in\mathcal Y$ is of the form
\[
(m,p)=\bigl(\E[y],\,\E[\mathbf 1_A(y)]\bigr)
\quad\text{for some }y\in\Sel(Y),
\]
this becomes
\[
U(\kappa)
=
\inf_{\lambda\in\R}\sup_{y\in\Sel(Y)}
\Bigl\{\E[\mathbf 1_A(y)]+\lambda(\E[y]-\kappa)\Bigr\}
=
\inf_{\lambda\in\R}\Bigl\{\E[\Psi(\lambda,\cdot)]-\lambda\kappa\Bigr\}.
\]
This proves the first identity in \eqref{eq:dual-envelopes-kappa}.

The formula for $L(\kappa)$ is obtained in the same way by replacing the
supremum problem for $\PP(y\in A)$ with the corresponding infimum problem.
Finally, since $U$ is the pointwise infimum of affine functions of $\kappa$,
it is concave, and since $L$ is the pointwise supremum of affine functions of
$\kappa$, it is convex.
\end{proof}

The dual envelopes arise from maximizing/minimizing the Lagrangian $x\mapsto \1\{x\in A\}+\lambda x$ over $x\in Y(\omega)$ and then interchanging expectation and extremum. The corresponding pointwise optimizers have a simple threshold structure.

\begin{theorem}[Threshold optimizer and calibration]\label{thm:threshold-calibration}
Assume $A$ is closed. For each $\lambda\in\R$, define the argmax multifunction
\[
Y_\lambda(\omega)
:=
\arg\max_{x\in Y(\omega)}\bigl(\1_{\{x\in A\}}+\lambda x\bigr).
\]
Then $Y_\lambda$ has nonempty closed values and admits measurable selections. Moreover, its graph is measurable.

If $\lambda>0$, then
\[
Y_\lambda(\omega)=
\begin{cases}
\{y_U(\omega)\}, & Y(\omega)\cap A=\emptyset,\\[4pt]
\{y_U(\omega)\}, & Y(\omega)\cap A\neq\emptyset,\ \lambda\Delta^+(\omega)>1,\\[4pt]
\{a^+(\omega)\}, & Y(\omega)\cap A\neq\emptyset,\ \lambda\Delta^+(\omega)<1,\\[4pt]
\{a^+(\omega),\,y_U(\omega)\}, & Y(\omega)\cap A\neq\emptyset,\ \lambda\Delta^+(\omega)=1.
\end{cases}
\]

Hence, outside the tie set
\[
T_\lambda^+
:=
\{\omega:Y(\omega)\cap A\neq\emptyset,\ \lambda\Delta^+(\omega)=1\},
\]
every measurable maximizer satisfies
\begin{equation}\label{eq:threshold-form-plus}
y_\lambda(\omega)
=
\begin{cases}
y_U(\omega), & Y(\omega)\cap A=\emptyset,\\[4pt]
y_U(\omega)-\min\{\Delta^+(\omega),\,1/\lambda\},
  & Y(\omega)\cap A\neq\emptyset.
\end{cases}
\end{equation}
Note that on $\{Y\cap A\neq\emptyset\}\setminus T_\lambda^+$, the
formula evaluates to $a^+(\omega)$ when
$\lambda\Delta^+(\omega)<1$ and to $y_U(\omega)$ when
$\lambda\Delta^+(\omega)>1$ (since
$\min\{\Delta^+,1/\lambda\}=\Delta^+=y_U-a^+$ in the first case
and $=1/\lambda$ in the second, with $1/\lambda<\Delta^+\le y_U-y_L$).
When $Y(\omega)\cap A=\emptyset$, the unique maximizer is $y_U(\omega)$
as stated.
If $\lambda<0$, then
\[
Y_\lambda(\omega)=
\begin{cases}
\{y_L(\omega)\}, & Y(\omega)\cap A=\emptyset,\\[4pt]
\{y_L(\omega)\}, & Y(\omega)\cap A\neq\emptyset,\ (-\lambda)\Delta^-(\omega)>1,\\[4pt]
\{a^-(\omega)\}, & Y(\omega)\cap A\neq\emptyset,\ (-\lambda)\Delta^-(\omega)<1,\\[4pt]
\{a^-(\omega),\,y_L(\omega)\}, & Y(\omega)\cap A\neq\emptyset,\ (-\lambda)\Delta^-(\omega)=1.
\end{cases}
\]
Hence, outside the tie set
\[
T_\lambda^-
:=
\{\omega:Y(\omega)\cap A\neq\emptyset,\ (-\lambda)\Delta^-(\omega)=1\},
\]
every measurable maximizer satisfies

\begin{equation}\label{eq:threshold-form-minus}
y_\lambda(\omega)
=
\begin{cases}
y_L(\omega), & Y(\omega)\cap A=\emptyset,\\[4pt]
y_L(\omega)+\min\{\Delta^-(\omega),\,1/(-\lambda)\},
  & Y(\omega)\cap A\neq\emptyset.
\end{cases}
\end{equation}

If $\lambda=0$, then
\[
Y_0(\omega)=
\begin{cases}
Y(\omega)\cap A, & Y(\omega)\cap A\neq\emptyset,\\
Y(\omega), & Y(\omega)\cap A=\emptyset.
\end{cases}
\]
In particular, $Y_0$ is again a measurable closed-valued multifunction.

Moreover, for every $\kappa\in[\E y_L,\E y_U]$, there exist $\lambda^\star\in\R$ and
$y^\star\in\Sel(Y_{\lambda^\star})$ such that
\[
\E[y^\star]=\kappa,
\qquad
\PP(y^\star\in A)=U(\kappa).
\]
\end{theorem}

\begin{proof}
Fix $\lambda\in\R$ and define
\[
f_\lambda(\omega,x):=\1_A(x)+\lambda x,
\qquad x\in Y(\omega).
\]

We first identify the argmax multifunction $Y_\lambda$ pointwise.

\medskip\noindent
\emph{Case 1: $\lambda>0$.}
If $Y(\omega)\cap A=\emptyset$, then $f_\lambda(\omega,x)=\lambda x$ on $Y(\omega)$,
hence the unique maximizer is $y_U(\omega)$.

Now assume $Y(\omega)\cap A\neq\emptyset$. Since $A$ is closed and $Y(\omega)$ is compact,
\[
a^+(\omega):=\sup\bigl(Y(\omega)\cap A\bigr)\in Y(\omega)\cap A.
\]
On $Y(\omega)\cap A$, the map $x\mapsto 1+\lambda x$ is increasing, so its maximum is
attained at $a^+(\omega)$ and equals $1+\lambda a^+(\omega)$. On $Y(\omega)\setminus A$,
the map $x\mapsto \lambda x$ is increasing, so its maximum is attained at $y_U(\omega)$
and equals $\lambda y_U(\omega)$. Thus the maximizers are determined by comparing
\[
\lambda y_U(\omega)
\quad\text{and}\quad
1+\lambda a^+(\omega).
\]
Since
\[
\lambda y_U(\omega)\ge 1+\lambda a^+(\omega)
\iff
\lambda\bigl(y_U(\omega)-a^+(\omega)\bigr)\ge 1
\iff
\lambda\Delta^+(\omega)\ge 1,
\]
we obtain exactly the stated description of $Y_\lambda(\omega)$ for $\lambda>0$.
Outside the tie set $T_\lambda^+$, this is equivalent to
\eqref{eq:threshold-form-plus}.

\medskip\noindent
\emph{Case 2: $\lambda<0$.}
If $Y(\omega)\cap A=\emptyset$, then $f_\lambda(\omega,x)=\lambda x$ on $Y(\omega)$,
and since $\lambda<0$ the unique maximizer is $y_L(\omega)$.

Now assume $Y(\omega)\cap A\neq\emptyset$. Since $A$ is closed and $Y(\omega)$ is compact,
\[
a^-(\omega):=\inf\bigl(Y(\omega)\cap A\bigr)\in Y(\omega)\cap A.
\]
On $Y(\omega)\cap A$, the map $x\mapsto 1+\lambda x$ is decreasing, so its maximum is
attained at $a^-(\omega)$ and equals $1+\lambda a^-(\omega)$. On $Y(\omega)\setminus A$,
the map $x\mapsto \lambda x$ is decreasing, so its maximum is attained at $y_L(\omega)$
and equals $\lambda y_L(\omega)$. Hence the maximizers are determined by
\[
\lambda y_L(\omega)\ge 1+\lambda a^-(\omega)
\iff
(-\lambda)\bigl(a^-(\omega)-y_L(\omega)\bigr)\ge 1
\iff
(-\lambda)\Delta^-(\omega)\ge 1.
\]
This gives the stated description of $Y_\lambda(\omega)$ for $\lambda<0$.
Outside the tie set $T_\lambda^-$, this is equivalent to
\eqref{eq:threshold-form-minus}.

\medskip\noindent
\emph{Case 3: $\lambda=0$.}
Then $f_0(\omega,x)=\1_A(x)$, so the maximizers are exactly
\[
Y_0(\omega)=
\begin{cases}
Y(\omega)\cap A, & Y(\omega)\cap A\neq\emptyset,\\
Y(\omega), & Y(\omega)\cap A=\emptyset.
\end{cases}
\]

In all three cases, $Y_\lambda(\omega)$ is a nonempty closed subset of the compact interval
$Y(\omega)$. The preceding pointwise characterization shows that $Y_\lambda$ is built from
the measurable objects $y_L$, $y_U$, $a^+$, $a^-$, $\Delta^+$, and $\Delta^-$ by measurable
case distinctions; hence $\operatorname{graph}(Y_\lambda)\in\mathcal A\otimes\mathcal B(\R)$.
Therefore the Kuratowski--Ryll-Nardzewski theorem yields $\Sel(Y_\lambda)\neq\emptyset$.

\medskip
Now let $y_\lambda\in\Sel(Y_\lambda)$. By construction,
\[
f_\lambda(\omega,y_\lambda(\omega))
=
\sup_{x\in Y(\omega)} f_\lambda(\omega,x)
=
\Psi(\lambda,\omega)
\qquad\text{a.s.}
\]
Therefore
\begin{equation}\label{eq:maximizer-lagrangian}
\E\big[\1_A(y_\lambda)+\lambda y_\lambda\big]
=
\E[\Psi(\lambda,\cdot)].
\end{equation}
In particular, every measurable selection of $Y_\lambda$ attains the pointwise-maximized
Lagrangian.

\medskip
Fix $\kappa\in[\E y_L,\E y_U]$ and define
\[
g(\lambda):=\E[\Psi(\lambda,\cdot)]-\lambda\kappa.
\]
By Theorem~\ref{thm:dual},
\[
U(\kappa)=\inf_{\lambda\in\R} g(\lambda).
\]
Since for each $\omega$ the map $\lambda\mapsto \Psi(\lambda,\omega)$ is the supremum of
affine functions of $\lambda$, $g$ is convex and lower semicontinuous.

If $\kappa=\E y_U$, then $y_U\in\Sel(Y\mid\kappa)$ and
\[
U(\kappa)=\PP(y_U\in A),
\]
so $y^\star:=y_U$ is optimal. Likewise, if $\kappa=\E y_L$, then $y^\star:=y_L$ is optimal.
Hence it remains to treat the interior case
\[
\kappa\in(\E y_L,\E y_U).
\]
For $\lambda\ge 0$,
\[
\Psi(\lambda,\omega)\ge \lambda y_U(\omega),
\]
so
\[
g(\lambda)\ge \lambda(\E y_U-\kappa)\xrightarrow[\lambda\to+\infty]{}+\infty.
\]
For $\lambda<0$,
\[
\Psi(\lambda,\omega)\ge \lambda y_L(\omega),
\]
so
\[
g(\lambda)\ge \lambda(\E y_L-\kappa)\xrightarrow[\lambda\to-\infty]{}+\infty.
\]
Thus $g$ is coercive and attains its minimum at some $\lambda^\star\in\R$.

\medskip
We next show that one can calibrate the mean at $\lambda^\star$.

For $\lambda>0$, define the lower and upper tie-breaking selections
\[
y_\lambda^-(\omega):=
\begin{cases}
a^+(\omega), & \omega\in T_\lambda^+,\\
\text{the unique maximizer in }Y_\lambda(\omega), & \omega\notin T_\lambda^+,
\end{cases}
\]
and
\[
y_\lambda^+(\omega):=
\begin{cases}
y_U(\omega), & \omega\in T_\lambda^+,\\
\text{the unique maximizer in }Y_\lambda(\omega), & \omega\notin T_\lambda^+.
\end{cases}
\]
For $\lambda<0$, define analogously
\[
y_\lambda^-(\omega):=
\begin{cases}
y_L(\omega), & \omega\in T_\lambda^-,\\
\text{the unique maximizer in }Y_\lambda(\omega), & \omega\notin T_\lambda^-,
\end{cases}
\]
and
\[
y_\lambda^+(\omega):=
\begin{cases}
a^-(\omega), & \omega\in T_\lambda^-,\\
\text{the unique maximizer in }Y_\lambda(\omega), & \omega\notin T_\lambda^-.
\end{cases}
\]
Finally, for $\lambda=0$, choose measurable selections $y_0^-,y_0^+\in\Sel(Y_0)$ such that
\[
y_0^-(\omega)=
\begin{cases}
\inf(Y(\omega)\cap A), & Y(\omega)\cap A\neq\emptyset,\\
y_L(\omega), & Y(\omega)\cap A=\emptyset,
\end{cases}
\]
and
\[
y_0^+(\omega)=
\begin{cases}
\sup(Y(\omega)\cap A), & Y(\omega)\cap A\neq\emptyset,\\
y_U(\omega), & Y(\omega)\cap A=\emptyset.
\end{cases}
\]

Set
\[
m_-(\lambda):=\E[y_\lambda^-],
\qquad
m_+(\lambda):=\E[y_\lambda^+].
\]
By construction,
\[
m_-(\lambda)\le \E[y]\le m_+(\lambda)
\qquad\text{for every }y\in\Sel(Y_\lambda).
\]

We claim that
\begin{equation}\label{eq:onesided-derivatives}
g'_-(\lambda)=m_-(\lambda)-\kappa,
\qquad
g'_+(\lambda)=m_+(\lambda)-\kappa.
\end{equation}
Indeed, for each fixed $\omega$, the one-sided derivatives of
$\lambda\mapsto\Psi(\lambda,\omega)$ are exactly the smallest and largest slopes
of the active affine pieces, namely
\[
D_-\Psi(\lambda,\omega)=y_\lambda^-(\omega),
\qquad
D_+\Psi(\lambda,\omega)=y_\lambda^+(\omega).
\]
Moreover,
\[
|D_\pm\Psi(\lambda,\omega)|\le |y_L(\omega)|+|y_U(\omega)|
\qquad\text{a.s.},
\]
and the right-hand side is integrable. Hence dominated convergence yields
\[
\frac{d^-}{d\lambda}\E[\Psi(\lambda,\cdot)]
=
\E[D_-\Psi(\lambda,\cdot)]
=
m_-(\lambda),
\]
and
\[
\frac{d^+}{d\lambda}\E[\Psi(\lambda,\cdot)]
=
\E[D_+\Psi(\lambda,\cdot)]
=
m_+(\lambda).
\]
Subtracting the derivative of $-\lambda\kappa$ gives \eqref{eq:onesided-derivatives}.

Since $\lambda^\star$ minimizes the convex function $g$, we have
\[
g'_-(\lambda^\star)\le 0\le g'_+(\lambda^\star).
\]
Using \eqref{eq:onesided-derivatives},
\[
m_-(\lambda^\star)\le \kappa\le m_+(\lambda^\star).
\]

\medskip
We now construct a measurable $y^\star\in\Sel(Y_{\lambda^\star})$ with $\E[y^\star]=\kappa$.

If $\lambda^\star>0$, then every maximizer agrees outside $T_{\lambda^\star}^+$, while on
$T_{\lambda^\star}^+$ the only possible values are $a^+$ and $y_U$. On this tie set,
\[
y_U-a^+=\Delta^+=1/\lambda^\star
\qquad\text{a.s. on }T_{\lambda^\star}^+.
\]
Hence
\[
m_+(\lambda^\star)-m_-(\lambda^\star)
=
\frac{1}{\lambda^\star}\PP(T_{\lambda^\star}^+).
\]
If $\PP(T_{\lambda^\star}^+)=0$, then $m_-(\lambda^\star)=m_+(\lambda^\star)$ and
necessarily $\kappa=m_-(\lambda^\star)$; take $y^\star:=y_{\lambda^\star}^-$. If
$\PP(T_{\lambda^\star}^+)>0$, choose
\[
\theta:=\frac{\kappa-m_-(\lambda^\star)}
{m_+(\lambda^\star)-m_-(\lambda^\star)}\in[0,1].
\]
By non-atomicity there exists a measurable set $B\subseteq T_{\lambda^\star}^+$ with
\[
\PP(B)=\theta\,\PP(T_{\lambda^\star}^+).
\]
Define
\[
y^\star:=\1_B\,y_{\lambda^\star}^+ + \1_{B^c}\,y_{\lambda^\star}^-.
\]
Then $y^\star\in\Sel(Y_{\lambda^\star})$ and
\[
\E[y^\star]
=
m_-(\lambda^\star)
+\frac{1}{\lambda^\star}\PP(B)
=
m_-(\lambda^\star)
+\theta\bigl(m_+(\lambda^\star)-m_-(\lambda^\star)\bigr)
=
\kappa.
\]

If $\lambda^\star<0$, the same argument applies on $T_{\lambda^\star}^-$, where
\[
a^--y_L=\Delta^- = 1/(-\lambda^\star)
\qquad\text{a.s. on }T_{\lambda^\star}^-.
\]
Thus one again obtains a measurable $y^\star\in\Sel(Y_{\lambda^\star})$ with
$\E[y^\star]=\kappa$.

If $\lambda^\star=0$, then
\[
Y_0(\omega)=
\begin{cases}
Y(\omega)\cap A, & Y(\omega)\cap A\neq\emptyset,\\
Y(\omega), & Y(\omega)\cap A=\emptyset,
\end{cases}
\]
is a one-dimensional random closed set. Since
\[
m_-(0)\le \kappa\le m_+(0),
\]
the one-dimensional Aumann expectation formula yields a selection
$y^\star\in\Sel(Y_0)$ with $\E[y^\star]=\kappa$.

\medskip
Finally, since $y^\star\in\Sel(Y_{\lambda^\star})$, it is a measurable pointwise maximizer.
Hence by \eqref{eq:maximizer-lagrangian},
\[
\E\big[\1_A(y^\star)+\lambda^\star y^\star\big]
=
\E[\Psi(\lambda^\star,\cdot)].
\]
Using $\E[y^\star]=\kappa$, we get
\[
\PP(y^\star\in A)
=
\E[\Psi(\lambda^\star,\cdot)]-\lambda^\star\kappa
=
g(\lambda^\star)
=
\inf_{\lambda\in\R}g(\lambda)
=
U(\kappa),
\]
where the last equality is Theorem~\ref{thm:dual}. Thus $y^\star$ attains $U(\kappa)$.
\end{proof}
The next result identifies the optimal probability in closed form once the calibration threshold is determined.
\begin{corollary}[Extremal probability formula]\label{cor:IdProb}
Let $\lambda^\star\in\R$ and let
\[
y^\star\in\Sel(Y_{\lambda^\star})
\]
be any optimizer provided by Theorem~\ref{thm:threshold-calibration}, so that
\[
\E[y^\star]=\kappa,
\qquad
\PP(y^\star\in A)=U(\kappa).
\]
Set
\[
t^\star:=1/\lambda^\star \qquad \text{if } \lambda^\star>0,
\qquad
s^\star:=1/(-\lambda^\star) \qquad \text{if } \lambda^\star<0.
\]
Then
\begin{equation}\label{eq:Id-prob}
\PP(y^\star\in A)=
\begin{cases}
\PP\!\big(\Delta^+< t^\star\big)+\eta^+\,\PP\!\big(\Delta^+=t^\star\big),
& \lambda^\star>0,\\[4pt]
\PP\!\big(\Delta^-< s^\star\big)+\eta^-\,\PP\!\big(\Delta^-=s^\star\big),
& \lambda^\star<0,\\[4pt]
\PP\big(Y\cap A\neq\varnothing\big), & \lambda^\star=0,
\end{cases}
\end{equation}
for some $\eta^+,\eta^-\in[0,1]$ determined by the tie-breaking on the boundary
sets. Equivalently, the optimal probability is obtained by full inclusion of the
strict-threshold region and possible randomization on the tie set.

Moreover, the threshold parameters are pinned down by the mean-calibration conditions
\begin{equation}\label{eq:calibration-plus}
\kappa\in
\Big[
\E\!\big[y_U-\min\{\Delta^+,t^\star\}\big],\,
\E\!\big[y_U-\min\{\Delta^+,t^\star\}\big]
+\frac{1}{\lambda^\star}\PP(\Delta^+=t^\star)
\Big]
\qquad (\lambda^\star>0),
\end{equation}
and
\begin{equation}\label{eq:calibration-minus}
\kappa\in
\Big[
\E\!\big[y_L+\min\{\Delta^-,s^\star\}\big],\,
\E\!\big[y_L+\min\{\Delta^-,s^\star\}\big]
+\frac{1}{-\lambda^\star}\PP(\Delta^-=s^\star)
\Big]
\qquad (\lambda^\star<0).
\end{equation}
In particular, if
\[
\PP(\Delta^+=t^\star)=0 \quad \text{or} \quad \PP(\Delta^-=s^\star)=0,
\]
then no tie-breaking is needed and the calibration conditions reduce to the equalities
\[
\EE\!\big[y_U-\min\{\Delta^+,t^\star\}\big]=\kappa
\qquad (\lambda^\star>0),
\]
and
\[
\EE\!\big[y_L+\min\{\Delta^-,s^\star\}\big]=\kappa
\qquad (\lambda^\star<0).
\]
\end{corollary}

\begin{proof}
Suppose first that $\lambda^\star>0$. By Theorem~\ref{thm:threshold-calibration}, the tie set is
\[
T_{\lambda^\star}^+
=
\{\omega:Y(\omega)\cap A\neq\emptyset,\ \lambda^\star\Delta^+(\omega)=1\}
=
\{\omega:\Delta^+(\omega)=t^\star\},
\]
and outside this set every maximizer is uniquely given by
\[
y^\star(\omega)=y_U(\omega)-\min\{\Delta^+(\omega),t^\star\}.
\]
Hence:
\[
\Delta^+(\omega)<t^\star \Longrightarrow y^\star(\omega)=a^+(\omega)\in A,
\]
\[
\Delta^+(\omega)>t^\star \Longrightarrow y^\star(\omega)=y_U(\omega)\notin A,
\]
while on $\{\Delta^+=t^\star\}$ both $a^+(\omega)\in A$ and $y_U(\omega)\notin A$
are pointwise maximizers. Therefore
\[
\PP(y^\star\in A)
=
\PP(\Delta^+<t^\star)+\eta^+\PP(\Delta^+=t^\star)
\]
for some $\eta^+\in[0,1]$, depending on the measurable tie-breaking on the boundary set.

For the mean, the minimal choice among maximizers is
\[
y_\star^-(\omega):=y_U(\omega)-\min\{\Delta^+(\omega),t^\star\},
\]
which selects $a^+$ on the tie set, while the maximal choice selects $y_U$ on the tie set.
Since on $\{\Delta^+=t^\star\}$,
\[
y_U(\omega)-a^+(\omega)=\Delta^+(\omega)=t^\star=\frac{1}{\lambda^\star},
\]
the set of attainable means across measurable tie-breakings is exactly
\[
\Big[
\E\!\big[y_U-\min\{\Delta^+,t^\star\}\big],\,
\E\!\big[y_U-\min\{\Delta^+,t^\star\}\big]
+\frac{1}{\lambda^\star}\PP(\Delta^+=t^\star)
\Big].
\]
Since $y^\star$ is chosen so that $\E[y^\star]=\kappa$, this proves
\eqref{eq:calibration-plus}.

The case $\lambda^\star<0$ is analogous. Now the tie set is
\[
T_{\lambda^\star}^-
=
\{\omega:Y(\omega)\cap A\neq\emptyset,\ (-\lambda^\star)\Delta^-(\omega)=1\}
=
\{\omega:\Delta^-=s^\star\},
\]
and outside this set every maximizer is uniquely given by
\[
y^\star(\omega)=y_L(\omega)+\min\{\Delta^-(\omega),s^\star\}.
\]
Thus
\[
\PP(y^\star\in A)
=
\PP(\Delta^-<s^\star)+\eta^-\PP(\Delta^-=s^\star)
\]
for some $\eta^-\in[0,1]$. Moreover, on $\{\Delta^-=s^\star\}$,
\[
a^-(\omega)-y_L(\omega)=\Delta^-(\omega)=s^\star=\frac{1}{-\lambda^\star},
\]
so the attainable means across measurable tie-breakings form exactly the interval
\[
\Big[
\E\!\big[y_L+\min\{\Delta^-,s^\star\}\big],\,
\E\!\big[y_L+\min\{\Delta^-,s^\star\}\big]
+\frac{1}{-\lambda^\star}\PP(\Delta^-=s^\star)
\Big],
\]
which gives \eqref{eq:calibration-minus}.

Finally, if $\lambda^\star=0$, then by Theorem~\ref{thm:threshold-calibration},
\[
Y_0(\omega)=
\begin{cases}
Y(\omega)\cap A, & Y(\omega)\cap A\neq\emptyset,\\
Y(\omega), & Y(\omega)\cap A=\emptyset,
\end{cases}
\]
so every optimizer chooses a point in $A$ whenever possible. Hence
\[
\PP(y^\star\in A)=\PP(Y\cap A\neq\emptyset).
\]
This completes the proof.
\end{proof}
Finally, we record qualitative properties and a convenient specialization.

\begin{proposition}[Extremes and global consistency]\label{prop:extremes}
At the mean extremes,
\[
U(\EE y_U)=L(\EE y_U)=\PP\big(y_U\in A\big),\qquad
U(\EE y_L)=L(\EE y_L)=\PP\big(y_L\in A\big),
\]
and for all $\kappa\in[\EE y_L,\EE y_U]$ the global random-set bounds are preserved:
\[
\PP\!\big(Y\subseteq A\big)\ \le\ L(\kappa)\ \le\ U(\kappa)\ \le\ \PP\!\big(Y\cap A\neq\varnothing\big).
\]

\end{proposition}

\begin{proof}
When $\kappa=\E y_U$, the mean constraint forces $y=y_U$ almost surely, so
$U(\kappa)=L(\kappa)=\PP(y_U\in A)$; similarly for $\kappa=\E y_L$. The
global bounds follow because the mean restriction can only shrink the feasible
set of selections, and Proposition~\ref{prop:RS-unrestricted} applies to any
$y\in\Sel(Y\mid\kappa)$.
\end{proof}

\begin{remark}[Interval targets]\label{rem:interval-target}
When $A=[a,b]$,
\[
\Delta^+(\omega)=
\begin{cases}
\bigl(y_U(\omega)-b\bigr)_+, & y_L(\omega)\le b,\\
+\infty, & y_L(\omega)>b,
\end{cases}
\qquad
\Delta^-(\omega)=
\begin{cases}
\bigl(a-y_L(\omega)\bigr)_+, & y_U(\omega)\ge a,\\
+\infty, & y_U(\omega)<a.
\end{cases}
\]
Hence, in the high-mean regime $\lambda^\star>0$, the optimizer is described by
\eqref{eq:threshold-form-plus}: it switches from $y_U$ to the rightmost point of
$Y\cap A$ once the gap to the right edge $b$ is small enough. In the low-mean
regime $\lambda^\star<0$, the optimizer is described by
\eqref{eq:threshold-form-minus}: it switches from $y_L$ to the leftmost point of
$Y\cap A$ once the gap to the left edge $a$ is small enough. The corresponding
threshold $t^\star=1/\lambda^\star$ or $s^\star=1/(-\lambda^\star)$ is then pinned
down by the mean-calibration equations \eqref{eq:calibration-plus}--\eqref{eq:calibration-minus}.
\end{remark}
\subsection{Further extensions: higher moments and additional quantiles}\label{subsec:extensions}

We now indicate how the selection-based approach developed
above extends to higher moments and to quantile restrictions other than
the median. We keep the standing assumptions of Subsection~\ref{subsec:notation}
and of Section~\ref{sec:main}: $Y=[y_L,y_U]$ is a random interval with
$\PP(y_L\le y_U)=1$ and $\E(|y_L|+|y_U|)<\infty$ on a non-atomic probability
space $(\Omega,\mathcal A,\PP)$.

\subsubsection*{Higher moments}

We first consider restrictions on the $r$th moment of selections for some
fixed $r\in\R$.

\begin{definition}[Moment-restricted selection sets]\label{def:r-moment}
Let $r\in\R$ and $\mu_r\in\R$. The \emph{$r$th-moment restricted selection set} is
\[
\Sel_r(Y\mid\mu_r)
:=
\bigl\{y\in\Sel(Y):\E|y|^r<\infty,\ \E[y^r]=\mu_r\bigr\}.
\]
When no confusion arises we write $\Sel_r(\mu_r)$ for $\Sel_r(Y\mid\mu_r)$.
\end{definition}

We begin with a simple non-emptiness result in the case where the
pointwise power map $x\mapsto x^r$ preserves the interval structure of
$Y$.

\begin{proposition}[Non-emptiness under power transformations]\label{prop:r-moment-nonempty}
Suppose either
\begin{enumerate}[label=(\roman*),nosep]
\item $r$ is an odd integer, or
\item $r>0$ and $y_L\ge 0$ almost surely.
\end{enumerate}
Define the random set
\[
Y^{(r)}(\omega):=\{x^r:x\in Y(\omega)\},\qquad \omega\in\Omega.
\]
Then $Y^{(r)}$ is a random interval, and
\[
\E\bigl[Y^{(r)}\bigr]
=
\bigl[\E(y_L^r),\E(y_U^r)\bigr]
=
\overline{\bigl\{\E[y^r]:y\in\Sel(Y)\bigr\}}.
\]
In particular, $\Sel_r(Y\mid\mu_r)\neq\emptyset$ if and only if
$\mu_r\in[\E(y_L^r),\E(y_U^r)]$.
\end{proposition}

\begin{proof}
Under (i) the map $x\mapsto x^r$ is strictly increasing on $\R$, so
$Y^{(r)}(\omega)=[y_L(\omega)^r,y_U(\omega)^r]$ for all $\omega$.
Under (ii) the map is increasing on $[0,\infty)$ and $Y(\omega)\subseteq[0,\infty)$
almost surely, so again $Y^{(r)}(\omega)=[y_L(\omega)^r,y_U(\omega)^r]$.
In both cases $Y^{(r)}$ is a random interval with endpoints
$L_r:=y_L^r$ and $U_r:=y_U^r$. The Aumann expectation identity
\eqref{eq:Aumann-interval} applied to $Y^{(r)}$ yields
\[
\E\bigl[Y^{(r)}\bigr]=\bigl[\E(y_L^r),\E(y_U^r)\bigr].
\]

If $y\in\Sel(Y)$, then $y^r$ is a measurable selection of $Y^{(r)}$,
so $\E[y^r]\in\E[Y^{(r)}]$. Conversely, for any
$m\in\E[Y^{(r)}]$ there exists $z\in\Sel(Y^{(r)})$ with $\E[z]=m$.
Because $x\mapsto x^{1/r}$ is strictly increasing on the relevant
domain (all of $\R$ in case (i), $[0,\infty)$ in case (ii)), the random variable
$y:=z^{1/r}$ is a measurable selection of $Y$ satisfying $y^r=z$ and hence
$\E[y^r]=\E[z]=m$. Thus the closure of $\{\E[y^r]:y\in\Sel(Y)\}$ equals
$\E[Y^{(r)}]$, and the final statement follows from
Definition~\ref{def:r-moment}.
\end{proof}

We now describe the induced range of means under an $r$th-moment
restriction. Define, for $\mu_r\in[\E(y_L^r),\E(y_U^r)]$,
\[
\underline m(r,\mu_r)
:=
\inf_{y\in\Sel_r(Y\mid\mu_r)} \E[y],
\qquad
\overline m(r,\mu_r)
:=
\sup_{y\in\Sel_r(Y\mid\mu_r)} \E[y],
\]
with the convention that $\underline m(r,\mu_r)=+\infty$,
$\overline m(r,\mu_r)=-\infty$ if $\Sel_r(Y\mid\mu_r)=\emptyset$.

\begin{theorem}[Dual representation under an $r$th-moment restriction]\label{thm:r-moment-dual}
Let the assumptions of Proposition~\ref{prop:r-moment-nonempty} hold and
let $\mu_r\in[\E(y_L^r),\E(y_U^r)]$. For $\lambda\in\R$ and $\omega\in\Omega$
define
\[
\Psi_r(\lambda,\omega)
:=
\sup_{x\in Y(\omega)}\bigl(x+\lambda x^r\bigr),
\qquad
\Phi_r(\lambda,\omega)
:=
\inf_{x\in Y(\omega)}\bigl(x+\lambda x^r\bigr).
\]
Then
\begin{equation}\label{eq:r-moment-dual}
\overline m(r,\mu_r)
=
\inf_{\lambda\in\R}
\Bigl\{\E\bigl[\Psi_r(\lambda,\cdot)\bigr]-\lambda\mu_r\Bigr\},
\qquad
\underline m(r,\mu_r)
=
\sup_{\lambda\in\R}
\Bigl\{\E\bigl[\Phi_r(\lambda,\cdot)\bigr]-\lambda\mu_r\Bigr\}.
\end{equation}
In particular, the $r$th-moment restricted mean set
$\{\E[y]:y\in\Sel_r(Y\mid\mu_r)\}$ is a nonempty compact interval
$[\underline m(r,\mu_r),\overline m(r,\mu_r)]\subseteq[\E y_L,\E y_U]$.
\end{theorem}

\begin{proof}
We give the argument for $\overline m(r,\mu_r)$; the lower endpoint is
handled analogously.

Consider the optimization problem
\[
\sup\bigl\{\E[y]:y\in\Sel(Y),\ \E[y^r]=\mu_r\bigr\}.
\]
Introduce a Lagrange multiplier $\lambda\in\R$ for the moment constraint and
define the Lagrangian
\[
\mathcal L(y,\lambda)
:=
\E[y]+\lambda\bigl(\E[y^r]-\mu_r\bigr)
=
\E\bigl[y+\lambda y^r\bigr]-\lambda\mu_r,
\qquad y\in\Sel(Y).
\]
For fixed $\lambda$, maximizing $\mathcal L(\cdot,\lambda)$ over
$\Sel(Y)$ subject to $y(\omega)\in Y(\omega)$ is equivalent to
maximizing
\[
\E\bigl[f_\lambda(\omega,y(\omega))\bigr],
\qquad
f_\lambda(\omega,x):=x+\lambda x^r,
\]
over selections $y$ of $Y$. By an Artstein-type interchange theorem for
random closed sets (see \cite{Artstein1974,Molchanov2005}),
\[
\sup_{y\in\Sel(Y)} \E\bigl[f_\lambda(\omega,y(\omega))\bigr]
=
\E\Bigl[\sup_{x\in Y(\omega)} f_\lambda(\omega,x)\Bigr]
=
\E\bigl[\Psi_r(\lambda,\cdot)\bigr].
\]
Therefore
\[
\sup_{y\in\Sel(Y)}\mathcal L(y,\lambda)
=
\E\bigl[\Psi_r(\lambda,\cdot)\bigr]-\lambda\mu_r.
\]

On the other hand, for any $y\in\Sel_r(Y\mid\mu_r)$ one has
$\mathcal L(y,\lambda)=\E[y]$ for all $\lambda$, so
\[
\overline m(r,\mu_r)
=
\sup_{y\in\Sel_r(Y\mid\mu_r)}\E[y]
=
\sup_{y\in\Sel(Y)}\inf_{\lambda\in\R}\mathcal L(y,\lambda).
\]
The feasible set of pairs
\[
\bigl(\E[y],\E[y^r]\bigr),\qquad y\in\Sel(Y),
\]
is the compact convex subset of $\R^2$ given by the Aumann expectation
of the random closed set
\[
\Xi(\omega):=\{(x,x^r):x\in Y(\omega)\},
\]
cf.\ \cite{Artstein1974,Molchanov2005}. Standard convex duality arguments
(see, e.g., \cite{Rockafellar1970}) yield
\[
\sup_{y\in\Sel(Y)}\inf_{\lambda\in\R}\mathcal L(y,\lambda)
=
\inf_{\lambda\in\R}\sup_{y\in\Sel(Y)}\mathcal L(y,\lambda),
\]
so
\[
\overline m(r,\mu_r)
=
\inf_{\lambda\in\R}\sup_{y\in\Sel(Y)}\mathcal L(y,\lambda)
=
\inf_{\lambda\in\R}
\Bigl\{\E\bigl[\Psi_r(\lambda,\cdot)\bigr]-\lambda\mu_r\Bigr\},
\]
which is the first identity in \eqref{eq:r-moment-dual}. The expression for
$\underline m(r,\mu_r)$ follows analogously, by applying the same
argument to the minimization problem for $\E[y]$ with $\E[y^r]=\mu_r$.

Finally, because the set of feasible pairs $(\E[y],\E[y^r])$ is compact and
convex, the intersection with the affine constraint $\E[y^r]=\mu_r$ is either
empty or a compact interval in the first coordinate. Under the non-emptiness
condition of Proposition~\ref{prop:r-moment-nonempty}, the attainable mean
set $\{\E[y]:y\in\Sel_r(Y\mid\mu_r)\}$ is a nonempty compact interval contained
in $[\E y_L,\E y_U]$ by \eqref{eq:Aumann-interval}.
\end{proof}

\begin{remark}
For $r=1$ and $\mu_1=\kappa$ we have $\Sel_1(Y\mid\mu_1)=\Sel(Y\mid\kappa)$
and $\underline m(1,\mu_1)=\overline m(1,\mu_1)=\kappa$, recovering the
trivial point of the mean. For $r\neq 1$, the interval
$[\underline m(r,\mu_r),\overline m(r,\mu_r)]$ describes the shrinkage of
the benchmark range $[\E y_L,\E y_U]$ induced by the $r$th-moment restriction.
\end{remark}
\subsubsection*{Additional quantile restrictions}

We next consider quantile restrictions away from the median. Recall from
Proposition~\ref{prop:quantile_selection} that for each $\alpha\in(0,1)$ and each
$q_\alpha$ in the feasible range
\[
  q_\alpha \in \bigl[T_Y^{-1}(\alpha),\,C_Y^{-1}(\alpha)\bigr],
\]
there exists a selection $y\in\Sel(Y)$ such that $F_y^{-1}(\alpha)=q_\alpha$.

\begin{definition}[Quantile-restricted selection sets]\label{def:alpha-quantile}
Fix $\alpha\in(0,1)$ and $q_\alpha\in\R$. The \emph{$(\alpha,q_\alpha)$-quantile
restricted selection set} is
\[
  \Sel_\alpha(Y\mid q_\alpha)
  :=
  \bigl\{y\in\Sel(Y): F_y^{-1}(\alpha)=q_\alpha\bigr\},
\]
where $F_y^{-1}$ denotes the left-continuous quantile function of $y$.
\end{definition}

We are interested in the induced range of means
\[
  \Theta_E(\alpha,q_\alpha)
  :=
  \bigl\{\E[y]:\,y\in\Sel_\alpha(Y\mid q_\alpha)\bigr\}.
\]
We first record existence and basic bounds, and then show that
$\Theta_E(\alpha,q_\alpha)$ is a convex subset of $[\E y_L,\E y_U]$.

\begin{proposition}[Non-emptiness and basic bounds]\label{prop:alpha-quantile-mean}
Let $\alpha\in(0,1)$ and assume
\[
  q_\alpha\in\bigl[T_Y^{-1}(\alpha),\,C_Y^{-1}(\alpha)\bigr].
\]
Then $\Sel_\alpha(Y\mid q_\alpha)\neq\emptyset$, and every
$y\in\Sel_\alpha(Y\mid q_\alpha)$ satisfies
\[
  \E[y]\ \in\ \bigl[\E(y_L),\,\E(y_U)\bigr].
\]
In particular,
\[
  \emptyset\neq\Theta_E(\alpha,q_\alpha)\ \subseteq\ [\E y_L,\E y_U].
\]
\end{proposition}

\begin{proof}
By Proposition~\ref{prop:quantile_selection}, if
$q_\alpha\in[T_Y^{-1}(\alpha),C_Y^{-1}(\alpha)]$ then there exists a selection
$y\in\Sel(Y)$ such that $F_y^{-1}(\alpha)=q_\alpha$. This shows
$\Sel_\alpha(Y\mid q_\alpha)\neq\emptyset$.

For any $y\in\Sel(Y)$ we have $y_L\le y\le y_U$ almost surely, hence by
monotonicity of expectation,
\[
  \E(y_L)\ \le\ \E[y]\ \le\ \E(y_U).
\]
This applies in particular to $y\in\Sel_\alpha(Y\mid q_\alpha)$, so
\[
  \Theta_E(\alpha,q_\alpha)
  \subseteq
  [\E y_L,\E y_U].
\]
Together with non-emptiness this yields the claim.
\end{proof}

The quantile restriction is non-linear in $y$ pointwise, so $y^{(1)},y^{(2)}\in
\Sel_\alpha(Y\mid q_\alpha)$ does not imply that their convex combination
$\theta y^{(1)}+(1-\theta)y^{(2)}$ has the same $\alpha$-quantile. Convexity of
$\Theta_E(\alpha,q_\alpha)$ therefore has to be established indirectly, at the
level of distributions and by allowing exogenous randomization on a product
extension of the underlying probability space.

We first note that the class of distributions with a fixed
$\alpha$-quantile is convex.

\begin{lemma}[Convexity of the fixed-quantile constraint in law]\label{lem:quantile_convex_law}
Let $\alpha\in(0,1)$ and $q\in\R$. Let $F_1,F_2$ be distribution functions
on $\R$ such that
\[
  F_i^{-1}(\alpha)\ =\ q,\qquad i=1,2,
\]
where $F_i^{-1}(\alpha):=\inf\{t\in\R:F_i(t)\ge\alpha\}$ is the left-continuous
$\alpha$-quantile. For $\theta\in[0,1]$ define
\[
  F_\theta := \theta F_1 + (1-\theta)F_2.
\]
Then $F_\theta^{-1}(\alpha)=q$.
\end{lemma}

\begin{proof}
Write $q_\alpha(F)$ for the $\alpha$-quantile of a distribution function $F$.
The condition $q_\alpha(F_i)=q$ means
\[
  F_i(t)<\alpha\quad\text{for all }t<q,
  \qquad
  F_i(q)\ge\alpha,
  \qquad i=1,2.
\]

Fix $t<q$. Then $F_i(t)<\alpha$ for $i=1,2$, so
\[
  F_\theta(t)
  =
  \theta F_1(t)+(1-\theta)F_2(t)
  <
  \theta\alpha+(1-\theta)\alpha
  =
  \alpha.
\]
At $t=q$,
\[
  F_\theta(q)
  =
  \theta F_1(q)+(1-\theta)F_2(q)
  \ge
  \theta\alpha+(1-\theta)\alpha
  =
  \alpha.
\]
Thus $F_\theta(t)<\alpha$ for all $t<q$ and $F_\theta(q)\ge\alpha$. By the
definition of the left-continuous quantile,
\[
  F_\theta^{-1}(\alpha)
  =
  \inf\{t\in\R:F_\theta(t)\ge\alpha\}
  =
  q.
\]
\end{proof}

We now use this to show that $\Theta_E(\alpha,q_\alpha)$ is convex. The proof
formally works on a product extension that introduces independent randomization,
but because our  analysis only depends on the joint law of
$(Y,y)$, the attainable range of $\E[y]$ is unaffected by such an enlargement.

\begin{proposition}[Convexity of the mean range]\label{prop:alpha-quantile-convex}
Let $\alpha\in(0,1)$ and $q_\alpha$ be as in
Proposition~\ref{prop:alpha-quantile-mean}. Then $\Theta_E(\alpha,q_\alpha)$ is
a convex subset of $\R$, hence an interval (possibly degenerate) contained in
$[\E y_L,\E y_U]$. Equivalently, for any
$y^{(1)},y^{(2)}\in\Sel_\alpha(Y\mid q_\alpha)$ and any $\theta\in[0,1]$, there
exists $\tilde y\in\Sel_\alpha(Y\mid q_\alpha)$ such that
\[
  \E[\tilde y]
  =
  \theta\,\E[y^{(1)}]+(1-\theta)\,\E[y^{(2)}].
\]
\end{proposition}

\begin{proof}
Let $y^{(1)},y^{(2)}\in\Sel_\alpha(Y\mid q_\alpha)$ and
$\theta\in[0,1]$. Denote by $F_i$ the distribution function of $y^{(i)}$,
$i=1,2$. By definition of $\Sel_\alpha(Y\mid q_\alpha)$,
\[
  F_i^{-1}(\alpha)=q_\alpha,\qquad i=1,2.
\]

\smallskip\noindent
\emph{Step 1: product extension and randomization.}
Consider the product probability space
\[
  (\bar\Omega,\bar{\mathcal A},\bar\PP)
  :=
  (\Omega\times[0,1],\ \mathcal A\otimes\mathcal B([0,1]),\ \PP\otimes\lambda),
\]
where $\lambda$ denotes Lebesgue measure on $[0,1]$. Let $U$ be the canonical
$\mathrm{Unif}[0,1]$ random variable on $[0,1]$, independent of all
random objects defined on $(\Omega,\mathcal A,\PP)$. Extend the random interval
$Y$ to $\bar\Omega$ by
\[
  Y(\omega,u) := Y(\omega),\qquad (\omega,u)\in\bar\Omega.
\]
Each $y^{(i)}$ can be viewed as a random variable on $\bar\Omega$ by ignoring
the second coordinate:
\[
  y^{(i)}(\omega,u):=y^{(i)}(\omega),\qquad i=1,2.
\]
Then $y^{(i)}(\omega,u)\in Y(\omega,u)$ for all $(\omega,u)$, so
$y^{(i)}\in\Sel(Y)$ on the extended space as well.

\smallskip\noindent
\emph{Step 2: mixture selection.}
Define a randomized selection $\tilde y:\bar\Omega\to\R$ by
\[
  \tilde y(\omega,u)
  :=
  \begin{cases}
    y^{(1)}(\omega), & u\le\theta,\\[4pt]
    y^{(2)}(\omega), & u>\theta.
  \end{cases}
\]
For each $(\omega,u)$, $\tilde y(\omega,u)$ equals either $y^{(1)}(\omega)$ or
$y^{(2)}(\omega)$, so $\tilde y(\omega,u)\in Y(\omega,u)$ and
$\tilde y$ is a measurable selection of $Y$ on $\bar\Omega$.

By independence of $U$ and $(y^{(1)},y^{(2)})$, the distribution of $\tilde y$
under $\bar\PP$ is the mixture
\[
  \text{law}(\tilde y)
  =
  \theta\,\text{law}(y^{(1)})+(1-\theta)\,\text{law}(y^{(2)}),
\]
so its distribution function is
\[
  F_\theta(t)
  :=
  \theta F_1(t)+(1-\theta)F_2(t),\qquad t\in\R.
\]

\smallskip\noindent
\emph{Step 3: quantile and mean of the mixture.}
Since $F_i^{-1}(\alpha)=q_\alpha$ for $i=1,2$, Lemma~\ref{lem:quantile_convex_law}
applied with $q=q_\alpha$ yields
\[
  F_\theta^{-1}(\alpha)=q_\alpha.
\]
Thus $\tilde y$ satisfies the same $(\alpha,q_\alpha)$-quantile restriction:
$F_{\tilde y}^{-1}(\alpha)=q_\alpha$, i.e.\ $\tilde y\in\Sel_\alpha(Y\mid q_\alpha)$
when viewed on the extended space. Its mean under $\bar\PP$ is
\[
  \E_{\bar\PP}[\tilde y]
  =
  \theta\,\E_\PP[y^{(1)}]+(1-\theta)\,\E_\PP[y^{(2)}],
\]
because $U$ is independent and selects $y^{(1)}$ with probability $\theta$
and $y^{(2)}$ with probability $1-\theta$.

\smallskip\noindent
\emph{Step 4: interpretation on the original space.}
The pair $(Y,\tilde y)$ on $(\bar\Omega,\bar{\mathcal A},\bar\PP)$ has a
well-defined joint law supported on $\{(K,x)\in\mathcal K(\R)\times\R:x\in K\}$.
On any non-atomic probability space there exists a pair $(\hat Y,\hat y)$
with the same joint law, with $\hat Y$ distributed as $Y$ and $\hat y$ a
measurable selection of $\hat Y$. In particular, one may realize such a pair
on $(\Omega,\mathcal A,\PP)$ without changing the distributions or expectations
of the scalar functionals considered here. For this $\hat y$ we then have
\[
  \hat y\in\Sel_\alpha(Y\mid q_\alpha),
  \qquad
  \E[\hat y] = \E_{\bar\PP}[\tilde y]
  =
  \theta\,\E[y^{(1)}]+(1-\theta)\,\E[y^{(2)}].
\]

Since $y^{(1)},y^{(2)}$ and $\theta\in[0,1]$ were arbitrary, this shows that
$\Theta_E(\alpha,q_\alpha)$ is convex. Together with
Proposition~\ref{prop:alpha-quantile-mean}, which ensures that
$\Theta_E(\alpha,q_\alpha)$ is nonempty and contained in $[\E y_L,\E y_U]$, we
conclude that $\Theta_E(\alpha,q_\alpha)$ is an interval (possibly a singleton)
contained in $[\E y_L,\E y_U]$.
\end{proof}

It is convenient to summarize the induced mean range under a single
quantile restriction by its lower and upper endpoints.

\begin{definition}[quantile-restricted mean interval]\label{def:quantile-mean-interval}
For $\alpha\in(0,1)$ and $q_\alpha$ as above, define
\[
  \underline m(\alpha,q_\alpha)
  :=
  \inf \Theta_E(\alpha,q_\alpha),
  \qquad
  \overline m(\alpha,q_\alpha)
  :=
  \sup \Theta_E(\alpha,q_\alpha),
\]
with the convention that $\underline m(\alpha,q_\alpha)=+\infty$ and
$\overline m(\alpha,q_\alpha)=-\infty$ if $\Sel_\alpha(Y\mid q_\alpha)=\emptyset$.
Whenever $\Sel_\alpha(Y\mid q_\alpha)\neq\emptyset$, the closed interval
\[
  \bigl[\underline m(\alpha,q_\alpha),\,\overline m(\alpha,q_\alpha)\bigr]
  \subseteq
  [\E y_L,\E y_U]
\]
will be called the \emph{$(\alpha,q_\alpha)$-quantile-restricted mean interval}.
\end{definition}
By Proposition~\ref{prop:alpha-quantile-mean} and
Proposition~\ref{prop:alpha-quantile-convex}, $\Theta_E(\alpha,q_\alpha)$ is a
nonempty convex subset of $[\E y_L,\E y_U]$, so its infimum and supremum are
finite and satisfy $\underline m(\alpha,q_\alpha)\le\overline m(\alpha,q_\alpha)$.
The endpoints need not be attained by a selection in general, but for any
$\varepsilon>0$ there exist $y_\varepsilon^-,y_\varepsilon^+\in
\Sel_\alpha(Y\mid q_\alpha)$ such that
\[
  \E[y_\varepsilon^-]\le \underline m(\alpha,q_\alpha)+\varepsilon,
  \qquad
  \E[y_\varepsilon^+]\ge \overline m(\alpha,q_\alpha)-\varepsilon.
\]

\begin{remark}
For $\alpha=\tfrac12$ and a singleton median $q_{1/2}=m$, the set
$\Sel_{1/2}(Y\mid m)$ coincides with $\Sel(Y\mid m)$ from
Subsection~\ref{subsec:IdMean}, and the mean range
$\Theta_E(1/2,m)$ is exactly the interval
$[\E_{\min}(m),\E_{\max}(m)]$ characterized in
Proposition~\ref{prop:MedianRestrictedAumann}. For general
$\alpha\in(0,1)$, Propositions~\ref{prop:alpha-quantile-mean} and
\ref{prop:alpha-quantile-convex} show that a single quantile restriction
always leads to a well-defined convex set of means contained in
the Aumann range $[\E y_L,\E y_U]$. Extending the explicit quantile-integral
formulas of Proposition~\ref{prop:MedianRestrictedAumann} to arbitrary levels
$\alpha\neq\tfrac12$ would require a more refined analysis of the tails below
and above $q_\alpha$ and additional regularity assumptions; we do not pursue
this here.
\end{remark}

\subsection{Extension to random convex sets in \texorpdfstring{$\Rd$}{Rd}}

We now extend the one-dimensional results of the preceding subsections to random
compact convex sets $Y:\Omega\to\mathcal{K}_c(\Rd)$, under an approximation scheme
by random unions of dyadic cubes. The purpose of this subsection is not to claim
that every constrained problem in $\Rd$ admits an explicit closed form, but rather
to show that the analysis may be reduced, in an approximation sense, to
one-dimensional interval-selection problems.

Throughout this subsection, let
\[
Y:\Omega\to\mathcal K_c(\Rd)
\]
be a random compact convex set defined on the fixed complete non-atomic
probability space $(\Omega,\mathcal A,\PP)$, and assume
\begin{equation}\label{eq:convex-int-envelope}
\E\Big[\sup_{x\in Y(\omega)}\|x\|\Big]<\infty
\end{equation}
and
\begin{equation}\label{eq:full-dim}
\operatorname{int}(Y(\omega))\neq\emptyset
\qquad\text{for $\PP$-almost every }\omega,
\end{equation}
where $\operatorname{int}(\cdot)$ denotes the interior in $\Rd$.

\begin{definition}[Random cube and its selection set]\label{def:random-cube}
A \emph{random cube} in $\Rd$ is a random set of the form
\[
R(\omega)=\prod_{i=1}^d [y_{L,i}(\omega),y_{U,i}(\omega)],
\]
where each $[y_{L,i},y_{U,i}]$ is a random interval with
\[
y_{L,i}\le y_{U,i}\qquad\text{$\PP$-a.s.}
\]
and
\[
\E\big(|y_{L,i}|+|y_{U,i}|\big)<\infty.
\]
\end{definition}

\begin{proposition}[Product decomposition]\label{prop:cube-product}
Let
\[
R=\prod_{i=1}^d [y_{L,i},y_{U,i}]
\]
be a random cube. Then
\begin{equation}\label{eq:ext-cube}
\Sel(R)
=
\prod_{i=1}^d \Sel([y_{L,i},y_{U,i}]).
\end{equation}
Moreover, for any $\kappa=(\kappa_1,\dots,\kappa_d)\in\R^d$,
\begin{equation}\label{eq:ext-cube-constrained}
\Sel(R\mid \kappa)
=
\prod_{i=1}^d \Sel([y_{L,i},y_{U,i}]\mid \kappa_i).
\end{equation}
\end{proposition}

\begin{proof}
A measurable map $y=(y_1,\dots,y_d):\Omega\to\Rd$ belongs to $\Sel(R)$ if and only if
\[
y(\omega)\in R(\omega)\qquad\text{$\PP$-a.s.}
\]
By the Cartesian-product structure of $R(\omega)$, this is equivalent to
\[
y_i(\omega)\in [y_{L,i}(\omega),y_{U,i}(\omega)]
\qquad\text{for each }i=1,\dots,d,\ \text{$\PP$-a.s.}
\]
Hence \eqref{eq:ext-cube} holds.

For the mean-restricted set, the condition $\E[y]=\kappa$ is equivalent to
\[
\E[y_i]=\kappa_i,\qquad i=1,\dots,d.
\]
Combining this with \eqref{eq:ext-cube} yields \eqref{eq:ext-cube-constrained}.
\end{proof}

\begin{proposition}[Pasting over partitions]\label{prop:paste}
Let $\{R_j\}_{j=1}^N$ be random cubes with pairwise disjoint interiors, and define
\[
\mathcal R := \bigcup_{j=1}^N R_j.
\]
Then $y\in\Sel(\mathcal R)$ if and only if there exist a measurable partition
$\{A_j\}_{j=1}^N$ of $\Omega$ and selections
\[
y^{(j)}\in\Sel(R_j),\qquad j=1,\dots,N,
\]
such that
\begin{equation}\label{eq:ext-paste}
y(\omega)=\sum_{j=1}^N \1_{A_j}(\omega)\,y^{(j)}(\omega)
\qquad\text{$\PP$-a.s.}
\end{equation}
\end{proposition}

\begin{proof}
Suppose first that $y\in\Sel(\mathcal R)$. For each $j$, define
\[
A_j
:=
\{\omega:y(\omega)\in R_j(\omega)\}
\setminus
\bigcup_{k<j}\{\omega:y(\omega)\in R_k(\omega)\}.
\]
Since the interiors of the $R_j(\omega)$ are pairwise disjoint, the family
$\{A_j\}_{j=1}^N$ is a measurable partition of $\Omega$ up to null sets.
On each $A_j$, we have $y(\omega)\in R_j(\omega)$.

Choose arbitrary measurable selections $\tilde y^{(j)}\in\Sel(R_j)$.
Then define
\[
y^{(j)}(\omega)
:=
\1_{A_j}(\omega)\,y(\omega)+\1_{A_j^c}(\omega)\,\tilde y^{(j)}(\omega).
\]
Each $y^{(j)}$ is measurable and belongs to $\Sel(R_j)$.
By construction, \eqref{eq:ext-paste} holds.

Conversely, suppose \eqref{eq:ext-paste} holds for some measurable partition
$\{A_j\}_{j=1}^N$ and selections $y^{(j)}\in\Sel(R_j)$. Then for $\omega\in A_j$,
\[
y(\omega)=y^{(j)}(\omega)\in R_j(\omega)\subseteq \mathcal R(\omega).
\]
Hence $y\in\Sel(\mathcal R)$.
\end{proof}

\begin{definition}[Dyadic inner approximation]\label{def:dyadic}
For $n\in\mathbb N$ and $\alpha=(\alpha_1,\dots,\alpha_d)\in\mathbb Z^d$, let
\[
Q_\alpha^{(n)}
:=
\prod_{i=1}^d [\alpha_i 2^{-n},(\alpha_i+1)2^{-n}]
\]
denote the dyadic cube of side length $2^{-n}$ indexed by $\alpha$.
The \emph{$n$th dyadic inner approximation} of $Y$ is the random closed set
\[
\mathcal R^{(n)}(\omega)
:=
\bigcup\Big\{
Q_\alpha^{(n)}:\,
Q_\alpha^{(n)}\subseteq Y(\omega)
\Big\}.
\]
\end{definition}

Since each cube $Q_\alpha^{(n)}$ is deterministic and Borel, and since the event
$\{Q_\alpha^{(n)}\subseteq Y\}$ is measurable for measurable closed-valued random sets,
it follows that $\mathcal R^{(n)}$ is a measurable random closed set. Moreover,
\[
\mathcal R^{(n)}(\omega)\subseteq Y(\omega)
\qquad\text{for all }\omega.
\]
Hence
\[
\Sel(\mathcal R^{(n)})\subseteq \Sel(Y)
\qquad\text{for every }n.
\]

The geometric core of the approximation argument is the following pointwise lemma.

\begin{lemma}[Hausdorff approximation by dyadic inner
  unions]\label{lem:dyadic-Hausdorff}
Let $K\subset\Rd$ be a nonempty compact convex set
\textbf{with nonempty interior}, and for each $n\in\mathbb N$
define
\[
\mathcal R^{(n)}(K)
:=
\bigcup\Big\{Q_\alpha^{(n)}:\,Q_\alpha^{(n)}\subseteq K\Big\}.
\]
Then $\mathcal R^{(n)}(K)\neq\emptyset$ for all sufficiently large
$n$, $\mathcal R^{(n)}(K)\subseteq K$ for every $n$, and
\[
d_H\big(\mathcal R^{(n)}(K),K\big)\xrightarrow[n\to\infty]{}0.
\]
\end{lemma}
\begin{proof}
The inclusion $\mathcal R^{(n)}(K)\subseteq K$ is immediate from the definition.

It remains to prove Hausdorff convergence. Since $\mathcal R^{(n)}(K)\subseteq K$,
we only need to show that
\[
\sup_{x\in K}\operatorname{dist}\bigl(x,\mathcal R^{(n)}(K)\bigr)\to 0.
\]

Assume, to the contrary, that this fails. Then there exist $\varepsilon>0$, a sequence
$n_k\uparrow\infty$, and points $x_k\in K$ such that
\[
\operatorname{dist}\bigl(x_k,\mathcal R^{(n_k)}(K)\bigr)\ge \varepsilon
\qquad\text{for all }k.
\]
By compactness of $K$, after passing to a subsequence we may assume that
\[
x_k\to x\in K.
\]

Fix $\delta\in(0,\varepsilon/4)$. Since $K$ is convex and compact, for every
$x\in K$ and every $\delta>0$, the contracted point
\[
x_\delta:=(1-\delta)x+\delta x_0
\]
belongs to the relative interior of $K$, where $x_0$ is any fixed point in the relative
interior of $K$. Hence there exists $r_\delta>0$ such that
\[
B(x_\delta,r_\delta)\cap \operatorname{aff}(K)\subseteq K.
\]
Because the dyadic mesh tends to zero, for all sufficiently large $k$ there exists a dyadic
cube $Q_{\alpha_k}^{(n_k)}$ of side length $2^{-n_k}$ contained in
$B(x_\delta,r_\delta)\cap \operatorname{aff}(K)\subseteq K$ and satisfying
\[
\operatorname{dist}(x,Q_{\alpha_k}^{(n_k)})<2\delta.
\]
Hence, for all sufficiently large $k$,
\[
\operatorname{dist}\bigl(x,\mathcal R^{(n_k)}(K)\bigr)<2\delta.
\]
Using $x_k\to x$, we obtain
\[
\operatorname{dist}\bigl(x_k,\mathcal R^{(n_k)}(K)\bigr)
\le
\|x_k-x\|+\operatorname{dist}\bigl(x,\mathcal R^{(n_k)}(K)\bigr)
<
\delta+2\delta
<
\varepsilon
\]
for all sufficiently large $k$, a contradiction.

Therefore
\[
\sup_{x\in K}\operatorname{dist}\bigl(x,\mathcal R^{(n)}(K)\bigr)\to 0,
\]
and since $\mathcal R^{(n)}(K)\subseteq K$, this is equivalent to
\[
d_H\big(\mathcal R^{(n)}(K),K\big)\to 0.
\]
\end{proof}

\begin{proposition}[$L^1$-density of cube selections]\label{prop:density}
Let $Y:\Omega\to\mathcal K_c(\Rd)$ be a random compact convex set satisfying
\eqref{eq:convex-int-envelope} and~\eqref{eq:full-dim}, and let
$\mathcal R^{(n)}$ be its dyadic inner approximations from
Definition~\ref{def:dyadic}. Then
\begin{equation}\label{eq:ext-dense}
\Sel(Y)
=
\overline{\bigcup_{n\ge1}\Sel(\mathcal R^{(n)})}^{\,L^1(\PP;\Rd)}.
\end{equation}
\end{proposition}

\begin{proof}
Since $\mathcal R^{(n)}(\omega)\subseteq Y(\omega)$ for every $n$ and $\omega$, we have
\[
\bigcup_{n\ge1}\Sel(\mathcal R^{(n)})\subseteq \Sel(Y).
\]
Hence
\[
\overline{\bigcup_{n\ge1}\Sel(\mathcal R^{(n)})}^{\,L^1}
\subseteq \overline{\Sel(Y)}^{\,L^1}.
\]
Because $Y$ is closed-valued, $\Sel(Y)$ is closed in $L^1(\PP;\Rd)$ under
\eqref{eq:convex-int-envelope}; thus
\[
\overline{\bigcup_{n\ge1}\Sel(\mathcal R^{(n)})}^{\,L^1}
\subseteq \Sel(Y).
\]

It remains to prove the reverse inclusion. Let $y\in\Sel(Y)$ be arbitrary.
For each $n$, define the multifunction
\[
\Gamma_n(\omega)
:=
\arg\min_{z\in \mathcal R^{(n)}(\omega)} \|y(\omega)-z\|.
\]
Since $\mathcal R^{(n)}(\omega)$ is compact and nonempty whenever $Y(\omega)\neq\emptyset$,
the set $\Gamma_n(\omega)$ is nonempty and compact. Its graph is measurable, because
\[
\operatorname{graph}(\Gamma_n)
=
\Bigl\{
(\omega,z): z\in \mathcal R^{(n)}(\omega),\ 
\|y(\omega)-z\|=\operatorname{dist}(y(\omega),\mathcal R^{(n)}(\omega))
\Bigr\},
\]
and both $\operatorname{graph}(\mathcal R^{(n)})$ and the distance map are measurable.
Hence, by the Kuratowski--Ryll-Nardzewski selection theorem, there exists a measurable
selection
\[
z^{(n)}\in\Sel(\mathcal R^{(n)})
\]
such that
\[
\|y(\omega)-z^{(n)}(\omega)\|
=
\operatorname{dist}\bigl(y(\omega),\mathcal R^{(n)}(\omega)\bigr).
\]

By Lemma~\ref{lem:dyadic-Hausdorff}, for almost every $\omega$,
\[
d_H\bigl(\mathcal R^{(n)}(\omega),Y(\omega)\bigr)\to 0.
\]
Since $y(\omega)\in Y(\omega)$, this implies
\[
\|y(\omega)-z^{(n)}(\omega)\|
=
\operatorname{dist}\bigl(y(\omega),\mathcal R^{(n)}(\omega)\bigr)
\le
d_H\bigl(\mathcal R^{(n)}(\omega),Y(\omega)\bigr)
\to 0
\]
for almost every $\omega$.

Moreover,
\[
\|y(\omega)-z^{(n)}(\omega)\|
\le
\|y(\omega)\|+\|z^{(n)}(\omega)\|
\le
2\sup_{x\in Y(\omega)}\|x\|,
\]
because both $y(\omega)$ and $z^{(n)}(\omega)$ belong to $Y(\omega)$.
The right-hand side is integrable by \eqref{eq:convex-int-envelope}. Therefore
dominated convergence yields
\[
\|z^{(n)}-y\|_{L^1(\PP;\Rd)}\to 0.
\]
This proves that every $y\in\Sel(Y)$ belongs to the $L^1$-closure of
$\bigcup_{n\ge1}\Sel(\mathcal R^{(n)})$, and hence \eqref{eq:ext-dense} follows.
\end{proof}

The preceding proposition shows that arbitrary selections of random compact convex sets
can be approximated in $L^1$ by selections of random finite unions of cubes. The next
result explains how mean restrictions pass through this approximation scheme.

\begin{theorem}[Reduction to random intervals]\label{thm:reduction}
Let $Y:\Omega\to\mathcal K_c(\Rd)$ satisfy
\eqref{eq:convex-int-envelope} and~\eqref{eq:full-dim}, let
$\kappa\in\R^d$, and let $\{\mathcal R^{(n)}\}_{n\ge1}$ be the
dyadic inner approximations of Definition~\ref{def:dyadic}. Then
the following hold.

\begin{enumerate}[label=\textup{(\roman*)},nosep]
\item For each $n$, every selection problem on $\Sel(\mathcal R^{(n)})$ reduces to finitely
many one-dimensional interval-selection problems through
Propositions~\textup{\ref{prop:cube-product}} and~\textup{\ref{prop:paste}}.

\item If $y\in\Sel(Y\mid\kappa)$, then there exists a sequence
\[
z^{(n)}\in \Sel(\mathcal R^{(n)})
\]
such that
\[
z^{(n)}\to y\qquad\text{in }L^1(\PP;\Rd)
\]
and therefore
\[
\E[z^{(n)}]\to \kappa.
\]

\item Let $\varphi:\Sel(Y)\to\R$ be \(L^1\)-continuous. Then
\[
\sup_{y\in\Sel(Y\mid\kappa)}\varphi(y)
\le
\liminf_{n\to\infty}
\sup\Bigl\{\varphi(z): z\in\Sel(\mathcal R^{(n)}),\ \|\E[z]-\kappa\|\le \varepsilon_n\Bigr\}
\]
for every deterministic sequence $\varepsilon_n\downarrow0$ for which the feasible sets
on the right-hand side are nonempty for all large \(n\). The analogous statement holds
for the infimum.
\end{enumerate}
\end{theorem}

\begin{proof}
\emph{Part (i).}
Each $\mathcal R^{(n)}$ is a finite union of dyadic cubes. By
Proposition~\ref{prop:paste}, a selection of $\mathcal R^{(n)}$ is obtained by
pasting together selections of its cube components on a measurable partition of $\Omega$.
By Proposition~\ref{prop:cube-product}, each cube selection decomposes coordinate-wise into
one-dimensional interval selections. Thus every optimization over $\Sel(\mathcal R^{(n)})$
reduces to finitely many one-dimensional interval-selection problems.

\medskip\noindent
\emph{Part (ii).}
Let $y\in\Sel(Y\mid\kappa)$. By Proposition~\ref{prop:density}, there exists
$z^{(n)}\in\Sel(\mathcal R^{(n)})$ such that
\[
z^{(n)}\to y \qquad\text{in }L^1(\PP;\Rd).
\]
Since expectation is continuous on \(L^1(\PP;\Rd)\), it follows that
\[
\E[z^{(n)}]\to \E[y]=\kappa.
\]

\medskip\noindent
\emph{Part (iii).}
Fix \(y\in\Sel(Y\mid\kappa)\). By part (ii), there exists a sequence
\(z^{(n)}\in\Sel(\mathcal R^{(n)})\) such that
\[
z^{(n)}\to y \quad\text{in }L^1,
\qquad
\E[z^{(n)}]\to \kappa.
\]
Hence for every deterministic sequence \(\varepsilon_n\downarrow0\) with
\(\|\E[z^{(n)}]-\kappa\|\le \varepsilon_n\) eventually, we have
\[
z^{(n)}\in
\Bigl\{
z\in\Sel(\mathcal R^{(n)}):\ \|\E[z]-\kappa\|\le \varepsilon_n
\Bigr\}
\]
for all sufficiently large \(n\). Since \(\varphi\) is \(L^1\)-continuous,
\[
\varphi(z^{(n)})\to \varphi(y).
\]
Therefore
\[
\varphi(y)
\le
\liminf_{n\to\infty}
\sup\Bigl\{\varphi(z):z\in\Sel(\mathcal R^{(n)}),\ \|\E[z]-\kappa\|\le \varepsilon_n\Bigr\}.
\]
Taking the supremum over all \(y\in\Sel(Y\mid\kappa)\) yields the desired inequality.
The proof for the infimum is analogous.
\end{proof}

\begin{remark}[Full-dimensionality
  assumption]\label{rem:full-dim}
Condition~\eqref{eq:full-dim} ensures that the dyadic cubes
$Q_\alpha^{(n)}$ eventually fit inside $Y(\omega)$.  If
$Y(\omega)$ is lower-dimensional on a set of positive probability,
the $d$-dimensional dyadic inner approximation
$\mathcal R^{(n)}(\omega)$ may be empty for all $n$, and the
reduction scheme does not apply directly.

Two routes around this restriction are available.  First, if the
affine dimension $\dim\!\operatorname{aff}(Y(\omega))=k<d$ is
$\PP$-a.s.\ constant, one may apply an affine isometry to reduce
to a $k$-dimensional problem and use $k$-dimensional dyadic
cubes; the selection-set structure is preserved.  Second, one may
replace the inner approximation by an outer approximation
(enclosing $Y(\omega)$ in a union of cubes that shrink to
$Y(\omega)$ in Hausdorff distance) and obtain the analogous
density result from the opposite direction.  We do not pursue
these extensions here and restrict attention to the
full-dimensional case.
\end{remark}

\begin{remark}[What the reduction theorem gives]\label{rem:reduction-meaning}
The reduction theorem should be read as an approximation principle.

First, selections of $Y$ can be approximated in \(L^1\) by selections of dyadic inner
approximations \(\mathcal R^{(n)}\). Second, each \(\mathcal R^{(n)}\) decomposes into a
finite union of random cubes, and each cube decomposes coordinate-wise into interval
selection problems. Thus constrained optimization problems for general random compact
convex sets in \(\Rd\) may be approximated by finite families of one-dimensional
interval-selection problems.

What is obtained here is therefore a reduction scheme from the general convex-set
case to the interval case, rather than an exact closed-form representation of
\(\Sel(Y\mid\kappa)\) for arbitrary \(Y\subset\Rd\).
\end{remark}

\noindent\textbf{Credit authorship contribution statement: }

\textbf{Arie Beresteanu: } Writing - review \& editing, Writing - original draft, Methodology, Investigation, Conceptualization.

\textbf{Behrooz Moosavi Ramezanzadeh: } Writing - review \& editing, Writing - original draft, Methodology, Investigation, Conceptualization.

\textbf{ Data Availability:} No data are used in this research.

\textbf{ Funding:} This research is not funded through any grant. 

\textbf{ Conflicts of Interest:} The authors have no conflicts of interest to declare. 

\bibliographystyle{elsarticle-num}
\bibliography{ref}
\end{document}